\documentclass{amsart}
\date{}
\usepackage{amsmath}
\usepackage{amssymb}
\usepackage{amsfonts}
\usepackage{latexsym}
\usepackage{verbatim}
\numberwithin{equation}{section}
\newtheorem{lemma}[equation]{Lemma}
\newtheorem{proposition}[equation]{Proposition}
\newtheorem{theorem}[equation]{Theorem}

\newtheorem{conjecture}[equation]{Conjecture}
\newtheorem{prop}[equation]{Proposition}

\newtheorem{question}[equation]{Question}
\theoremstyle{remark}

\newtheorem{remark}[equation]{Remark}

\newtheorem{example}[equation]{Example}

\newcommand{\QED}{$\Box$}

\begin{document}

\title
[Periodic points] {Jointly periodic points in cellular automata:
computer explorations and conjectures }

\begin{abstract}
We develop  a rather 
elaborate computer program to investigate the 
jointly  periodic points of one-dimensional
cellular automata. The experimental results and mathematical context
lead to questions, conjectures and a contextual theorem. 
\end{abstract}

\author{Mike Boyle}
\address{Mike Boyle\\Department of Mathematics\\
        University of Maryland\\
        College Park, MD 20742-4015\\
        U.S.A.}
\email{mmb@math.umd.edu} 
\urladdr{www.math.umd.edu/$\sim$mmb}
\author{Bryant Lee} 
\address{Bryant Lee \\Department of Mathematics\\
        University of Maryland\\
        College Park, MD 20742-4015\\
        U.S.A.}
\email{blee3@umd.edu} 
\thanks{Partially supported by NSF Grants 0400493 and
NSFD.DMS0240049}



%below are the definitions for some subject classification numbers
%37B10 symbolic dynamics
%37B40 topological entropy
%37C40 smooth ergodic theory, invariant measures
%37C45 dimension theory of dynamical systems
%37C99 smooth dynamical systems, general theory
%37D35 thermodynamic formalism, variational principles, equilibrium states
%37D20 Uniformly hyperbolic systems (expanding, Anosov, Axiom A, etc.)
%37D25 Nonuniformly hyperbolic systems (Lyapunov exponents, Pesin
%theory,etc.)

\maketitle

\tableofcontents 

\section{Introduction and conjectures}

In this paper we consider the action of a
surjective one-dimensional cellular automaton  $f$  
on jointly  periodic points.  
Detailed definitions  are recalled below. 

This paper is primarily 
an experimental mathematics paper, based on 
data from a program written by the second-named author to 
explore such actions.  
The experimental results  and mathematical context 
lead us to questions and a conjecture 
on the growth rate of the jointly periodic points.   
The program itself is freely available at the website 
of the first-named author.

We approach our topic from the perspective of 
symbolic dynamics, which provides  some relevant 
tools and  results. However, almost   
all of this paper--in particular the questions and 
conjectures--can be well understood 
without symbolic dynamics.  
We do  spend time on  context, and even prove a theorem  
(Theorem \ref{bignu}), for two reasons. First, we believe 
that experimental mathematics should not be too segregated 
from the motivating and constraining mathematics. 
Second, workers on cellular automata have diverse backgrounds, 
not necessarily including symbolic dynamics. (Similarly, perhaps 
a technique or example unfamiliar to us  could  
resolve one of our questions.)

To express our questions and conjectures clearly, we must 
suffer some 
definitions. We let $\Sigma_N$ denote the set of doubly 
infinite sequences $x= \dots x_{-1}x_0x_1\dots $ such that 
each $x_i$ lies in a finite ``alphabet'' $\mathcal A$ of $N$ symbols;  
usually  $\mathcal A= \{0, 1, \dots , N-1\}$. A one-dimensional 
cellular automaton (c.a.) is a 
 function $f: \Sigma_N \to \Sigma_N$ for which there are integers 
$a\leq b$ and a function 
$F: \mathcal A^{b-a+1}\to \mathcal A$ 
($F$ is a ``local rule'' for $f$) such that 
 for all $i$, $(f(x))_i = F(x_{i+a}\cdots x_{i+b})$. 
The shift map $\sigma$ on a sequence  is defined by 
$(\sigma x)_i = x_{i+1}$.  
We let $S_N$ denote the shift map on $\Sigma_N$. 

For any map $S$, we 
let $P_k(S)$ denote the points of (not necessarily least) period $k$ 
of $S$, i.e. the points fixed by $S^k$, and let $\text{Per}(S)=
\cup_k P_k(S)$. 
Thus, $\text{Per}(S_N)$ is the set of ``spatially periodic'' points for a 
one-dimensional cellular automaton on $N$ symbols. 
The {\it jointly periodic} points of a cellular automaton map $f$ 
on $N$ symbols are the points in $\text{Per}(S)$ which are also  
periodic under $f$, that is, the points which are ``temporally
periodic'' as well as spatially periodic. (In the usual computer
screen display, 
this would mean vertically and well as horizontally periodic.)
There is by this time a lot of work addressing periodic and 
jointly periodic points for {\it linear} one-dimensional c.a.; we 
refer to \cite{CCG,CDN,Jen,LiTh,MOW,MST,Su} and their references. 
Also see \cite{Miles} regarding the structure of 
periodic points for these and more general algebraic maps 
in the setting of 
  \cite{kigroup,schm}.

A subset $E$ of $\Sigma_N$ is  {\it dense} if for every point 
$x$ in $\Sigma_N$ 
and every $k\in \mathbb N$ there exists $y$ in $E$ such that 
$x_i=y_i$ whenever $|i|\leq k$.  
We say  $E$ 
is $m$-dense if every word of length $m$ on symbols from 
the alphabet occurs in a point of $E$.  

We can now state our first conjecture. 

\begin{conjecture} \label{D} 
For every surjective one-dimensional cellular automaton, 
the jointly periodic points are dense. 
\end{conjecture} 

Conjecture \ref{D} is a known open question
\cite{Bl,BT,BKi},  justified by its clear relevance to
a dynamical systems approach to cellular automata. 
(Whether points which are temporally but not necessarily spatially 
periodic for a surjective c.a. must be dense is likewise unknown
\cite{Bl}.) 
That  this question, also open for higher dimensional c.a.,  
has not been answered reflects the difficulty of saying anything 
of a general nature about c.a., for which meaningful questions are 
often undecidable \cite{kari}.  

It is known that the jointly periodic points of a one-dimensional 
cellular automaton map $f$ are dense if $f$ is a closing map 
\cite{BKi}  or if $f$ is surjective with a point 
of equicontinuity \cite{BT}.   
  We justify our escalation of (\ref{D}) from question to conjecture 
by augmenting earlier results 
with some  experimental evidence. 
In particular: for every span 4 surjective one-dimensional cellular automaton 
on two symbols, the jointly periodic points are at least 13-dense 
(Proposition \ref{13dense}).

Now we turn to more quantitative questions. 
Letting for the moment $P$ denote
the number of points in 
$P_k(S_N)$ 
which are  periodic  
under $f$ as well as $S_N$ 
(i.e. $P = |\textnormal{Per}(f|P_k(S_N))|$), 
we set 
$\nu_k(f,S_N)=P^{1/k}$,  and then define 
\[
\nu (f,S_N)  = \limsup_k \nu_k(f,S_N) \ .
\] 
\begin{question}\label{Q1}
Is it true for every surjective one dimensional cellular automaton 
$f$ on $N$ symbols that 
$\nu (f,S_N) \geq \sqrt N$? 
\end{question}

\begin{question} \label{Q2}
Is it true for every surjective one dimensional cellular automaton 
$f$ on $N$ symbols that 
$\nu (f,S_N) > 1$? 
\end{question} 
We cannot answer Question \ref{Q2} even in the case 
that $f$ is a ``closing'' map 
and we know there is an abundance of jointly periodic points 
\cite{BKi}. 

\begin{conjecture} \label{C} 
There exists $N>1$ and a surjective cellular automaton $f$ 
on $N$ symbols such that 
$\nu (f,S_N) < N$.   
\end{conjecture} 
Conjecture  \ref{C} is a proclamation of ignorance. 
From the experimental data in our Tables, it seems 
perfectly clear that there will be many surjective 
c.a. $f$ with $\nu (f,S_N)< N$. However, we are unable to give a 
proof for any example.  With the additional assumption that 
the c.a. is linear, it is known that  Conjecture \ref{D} is true 
and the answer to Question \ref{Q1} is yes (Sec. \ref{mechanisms}). 

The relation of Questions (\ref{Q1}-\ref{C}) to Conjecture \ref{D} 
is the following: 
if a c.a. map $f$ on $N$ symbols does not have dense periodic 
points, then $\nu (f,S_N)<N$. 

Here is the organization of the sequel. 
In Section \ref{definitions}, 
we give detailed definitions and background. 
In Section \ref{mechanisms}, we establish some mechanisms by which 
one can prove lower bounds for $\nu(f,S_N)$ for some $f$. 
We also prove (Theorem \ref{bignu}) 
that no property of a surjective c.a. considered abstractly 
as a quotient map 
without iteration can  establish 
$\nu (f,S_N)<N$.  
We also support Question \ref{Q1}  with a random maps heuristic. 
(The potential analogy of c.a. and random maps 
was remarked earlier by Martin, Odlyzko 
and Wolfram \cite[p.252]{MOW} in their study of linear c.a.)   
A list of c.a. used for the computer explorations is given 
in Section \ref{maps}. 

Our computer program consists of three related subprograms: 
FDense, FPeriod and FProbPeriod. We use these respectively 
 in Sections \ref{secfdense} ,
\ref{FPeriod} and  \ref{secFProbPeriod}. 
FDense probes approximate density of jointly periodic points 
of a given shift period. 
FPeriod provides exact information on jointly periodic points 
of a given shift period. FProbPeriod provides information on 
jointly periodic points for a random sample from a given 
shift period, and thus provides some information at shift 
periods where the memory demands of FPeriod are too great for 
it to succeed.

In Sections 
\ref{secfdense} ,
\ref{FPeriod} and  \ref{secFProbPeriod}, 
we give more information on the algorithms  
and discuss the many tables of output data in the appendices.   
In the tables, decimal output data are approximated by truncation; 
e.g.,  1.429 becomes 1.42 rather than 1.43. 
Detailed information 
on the program is available, along with the program itself,  at the 
website of the first named author.

\section{Definitions and background}\label{definitions}

Let $\mathcal A=\{0, 1,\dots ,N-1 \}$,  a finite set of $N$ symbols, 
with the discrete
topology. Let   $\Sigma_N$ be the
product space
$\mathcal A^{\mathbb{Z}}$, with the product topology. We view a point
$x$ in $\Sigma_N$ as a doubly infinite sequence of
symbols from $\mathcal A$,  $x=\dots x_{-1}x_0x_1\dots $ .
The space $\Sigma_N$ is compact and metrizable; one metric
compatible with the topology is
$\text{dist}(x,y)= 1/(|n|+1)$ where
$|n|$ is the minimum nonnegative integer such that $x_n \neq y_n$.
A set $E$ is {\it dense} in $\Sigma_N$ in this topology if for every $k$ 
and every word $W$ in $\mathcal A^{2k+1}$ there exists $x$ in $E$ 
such that $x[-k,k]=W$. 

The shift map $\sigma$ sends a sequence $x$ to the sequence 
$\sigma x$ defined by $(\sigma x)_i=x_{i+1}$. The shift map 
defines a homeomorphism $S_N$ on $\Sigma_N$. 
The topological dynamical system $(\Sigma_N,S_N)$
is called the {\it full shift on $N$ symbols}, or more 
briefly the $N$-shift.  
%If $X$ is a closed shift-invariant subsystem of some full shift, 
%then $(X,\sigma|X)$ is a {\it subshift}. 
%IIf $X$ is exactly the 
%subsystem of points which never see some given finite set of 
%words, then $(X,\sigma|X)$ is a 
%{\it shift of finite type}. These subshifts are 
%useful and natural for the study of cellular 
%automata (see the remarks in \cite{BKi}). 
%For a lighter notation, 
%we may sometimes use a single letter like $S$ 
%to represent a subshift. 

%We will make use of the {\it entropy} of a subshift $S=(X,\sigma |X)$, 
%which  
%is $h(S)$, the growth rate of $S$-words:   
%$h(S) = \lim_k (1/k)\log (|\{x[1,k]: x\in X|\})$. 
%In particular, $h(S_N)=\log N$. 

A map $f: \Sigma_N\to \Sigma_N$ is  continuous and shift-commuting 
($f\sigma = \sigma f$) if and only $f$ is a {\it block code}, i.e. 
there exist integers $a,b$ 
and a function $F: \mathcal A^{b-a+1}\to \mathcal A$ such that 
$(f(x))_i = F(x[a,b])$ for integers 
$i$, for all $x\in \Sigma_N$.  Such a map $f$ is called a 
one-dimensional cellular automaton. 
There is a well known 
dichotomy for such maps $f$: either (i) $f$ is surjective and 
for some integer $M$ every point has at most $M$ preimages, 
or (ii) image points typically have uncountably many preimages, 
and $f$ is not surjective \cite{hedlund,Ki,LM}. 
In Case (i), almost all points have the same number of preimages; 
this number is the {\it degree} of $f$. 

We restrict our attention to surjective maps in this paper because 
we are interested in periodic points  of $f$, which  
must be contained in $\cap_{k>0} f^k\Sigma_N$, 
the eventual image of $f$.   
We separate our ignorance about periodic points from 
additional difficulties involving the passage to the 
eventual image \cite{Ma}. 

Polynomials can be used to define cellular automata; 
for example, if we refer to the c.a. $f$ 
defined on the $N$-shift by the polynomial 
$2x_{-1}+x_0(x_2)^3$, we mean that $f$ is defined by the block code
$(fx)_i = 2x_{i-1}+x_i(x_{i+2})^3$, where the arithmetic is 
interpreted modulo $N$.  The {\it span} 
of such a code is  1 plus the maximum difference of 
coordinates with nonzero coefficients; in this  
example, it is $1+2-(-1)=4$. The code is 
left permutative if for every $x$, permuting 
inputs to the leftmost variable, with inputs to other 
variables fixed,  permutes 
the outputs. Likewise there is the notion of right 
permutative. The previous example is left permutative
and it is not right permutative. 
When the number $N$ of symbols is prime, 
every c.a. map $f$ has such a polynomial representation
 \cite{hedlund}. 
(For general $N$, there is a representation by  
a product of polynomial  
representations over  finite fields  \cite{MOW}.) 

A block code on $S_N$ depending on coordinates $[0,j-1]$ 
can be described by a ``lookup code'', 
a word $W$ of length $N^j$ on alphabet
$\{0,...,j-1\}$ defined as follows. List the $N^j$ possible blocks 
of length $j$ in lexicographic order; then the $i$th 
symbol of $W$ is the output symbol under $f$ for the $i$th 
input block. For example, for the code $x_0 + x_1x_2$ on 
$S_2$, the input words in lexicographic order are 
000,001,010,011, 100,101,110,111 and the corresponding word 
$W$ is 0001 1110. 

For the $N$ shift, the number of coding rules of span at most 
$j$ is $N^{N^j}$. If $\textnormal{inj}(j,N)$ denotes the number of 
these which define injective (and thus surjective \cite{hedlund}) 
codes, then we still 
\cite{kimroush} 
see a superexponential growth rate in $j$,  
$\lim_j \log \log (\textnormal{inj}(j,N)) = \log N$, even though
surjective span $j$ maps become very sparse in the set of all span $j$
maps, as $j$ increases.   

A block code $f:\Sigma_N \to \Sigma_N$ is 
{\it right-closing} if it never collapses distinct left-asymptotic
points. 
 This means that  if $f (x)= f(x')$ and for some $I$
it holds that $x_i = x'_i$ for 
all $i$ in $(-\infty , I]$, then $x=x'$. 
Any right permutative map is right closing. 
The definition of {\it left closing} is given by replacing 
$(-\infty , I]$ with $[I,\infty )$. 
 The map $f$ is {\it closing} if it is either left or 
right closing. 
An endomorphism of a full shift $S_N$ is 
constant-to-one if and only if it is 
both right and left closing (i.e., it is biclosing).  
 A closing map is surjective. 
Closing maps are important in the coding theory of symbolic dynamics 
\cite{ashley,Ki,LM}. They also have a very natural description from the 
viewpoint of hyperbolic dynamics \cite{BrinStuck}: right closing maps are injective 
on unstable sets, left closing maps are injective on stable sets.

We now discuss some previous work involving periodic points and 
cellular automata. 
We let $P_n(S)$ denote the points of period 
$n$ of $S$, and $P^o_n(S)$ the points of least period $n$.
These finite sets are mapped into themselves by any c.a. map $f$; 
thus any periodic point of $S$ is at least preperiodic 
for $f$. For a preperiodic (possibly periodic) point $x$, 
the {\it preperiod} of $x$ is the least nonnegative integer $j$ such that 
$f^j(x)$ is periodic, and the {\it period }
of $x$ is its eventual period, the smallest positive integer $k$ 
such that $f^{m+k}(x)=f^m(x)$ for all large $m$. A point is jointly 
periodic if it is periodic under both $f$ and $S_N$. 

In the case $f$ is linear 
($f(x)+f(y)=f(x+y)$),  
Martin, Odlyzko and Wolfram  \cite{MOW}  
(see also  the further work in   
\cite{CCG,CDN,Jen,LiTh,MST,Su} and their references)
gave  an algebraic analysis  of $f$-periods and preperiods 
for points of a given shift period,
and also provided some numerical data. 
One key feature for linear $f$ is an easy observation: 
among the jointly periodic 
points of shift period $k$, there will be a  point (generally 
many points) whose least $f$-period will be an integer multiple 
of all the least $f$-periods of the jointly periodic points of 
shift period $k$. In contrast, 
a very special case of 
a powerful theorem of 
Ashley \cite{ashley} has the following 
statement:
 for any $K,N$ and 
any shift-commuting 
map $g$ from $\cup_{1\leq k \leq K}P_k(S_N)$ to itself, there 
will exist surjective c.a. on $N$ symbols whose restriction 
to $\cup_{1\leq k \leq K}P_k(S_N)$ equals $g$. 

The following remark is another indication of the difficulty of 
understanding joint periodicity of even injective c.a.  
For a map $T$, $\textnormal{Fix}(T)$ denotes $P_1(T)$, 
the set of fixed points of $T$.

\begin{remark} Given $N\geq 2$, let $S$ denote $S_N$, and   
suppose $\phi$ is an injective one-dimensional c.a. 
on $N$ symbols. Suppose $N$ is prime. 
Then there will exist some integer $m$ and some 
$\kappa > 0$ such that for all $k\in \mathbb N$, 
\[
\Big|
\textnormal{Fix} \big( (S^a\phi^b)^k \big)
\Big| 
= 
N^{(a+mb)k} 
= 
\Big|\textnormal{Fix} \big( (S^a(S^m)^b)^k \big)\Big|
= 
\Big|\textnormal{Fix} \big( (S^{a+mb})^k \big)\Big|
%\qquad \textnormal{for all {\it k} in } \mathbb Z
%\Big|P_k\big((S^m\phi)\big)\Big| 
%= 
%\Big|P_k\big((S^m)\big) \Big| \qquad \textnormal{for all {\it k} in } \mathbb Z
\] 
whenever $ |b/a| < \kappa $ 
(this follows from \cite[Theorem 2.17]{BKr1}). 
That is, for the two $\mathbb Z^2$ actions generated
respectively  by 
$S, \phi $ and $S,S^m$, the periodic point counts for 
actions by individual elements $(a,b)$ of $\mathbb Z^2$ 
are the same for 
all  $(a,b)$ in some open cone around the positive horizontal 
axis. Nevertheless, 
with the given $m$ fixed,  the sequence 
$(|P_k(\phi)|)$ can still vary tremendously with $\phi$.
(For a dramatic example in the setting of shifts of finite
type, see \cite[Example 10.1]{Na1}). 
\end{remark} 

Lastly, we note that the invariant $\nu$, defined in the introduction, 
 has an unusual robustness, 
as follows. 
\begin{remark} Fix $N$ and let $S=S_N$. 
Suppose $x\in \text{Per}(S_N)$ and $f$ 
is a c.a.  on $\Sigma_N$. Then $x$ is 
 in $\text{Per}(f)$ if and only if for some $i>0$,
$f^ix$ and $x$ are in the same $S$-orbit. It follows that 
for all integers $i,j,k$ with $k,i$ positive and $j$ 
nonnegative, we have 
$\nu_k(f,S)= \nu_k(f^iS^j, S)$, and thus 
 $\nu(f,S)= \nu(f^iS^j, S)$.  
%If also $k\in \mathbb N$, then 
%$\nu_{mk}(f,S)= \nu_{m}(f^iS^j, S^k)$, and therefore 
%$\nu(f,S^k)\leq  (\nu(f, S))^k$. 
\end{remark}

\section{Some mechanisms for periodicity} \label{mechanisms} 

Throughout this section $f$ denotes a c.a. map on $N$ symbols. 
In this section, we discuss 
four ways to prove $\nu (f,S_N)$ is large: 
\begin{enumerate} 
\item 
find a large shift fixed by $f$ (or more generally by 
a power of $f$)  
\item let $f$ be linear (i.e., a group endomorphism of $\Sigma_N$, 
where addition on the compact group $\Sigma_N$ is defined
coordinatewise mod $N$)  
%\cite{Ki0} 
\item use the algebra of a polynomial presenting $f$ 
\item 
find
equicontinuity points. 
\end{enumerate} 
After discussing these, 
we offer a random maps heuristic and a question. 

(1) We will exhibit the first mechanism in some generality. 
Two c.a. $f,g$ are isomorphic if there is an invertible c.a. 
$\phi$ such that $f=\phi g \phi^{-1}$ (where e.g. 
$\phi g$ is the composition,  $(\phi g)(x)=\phi (g(x))$).  
The c.a. $f,g$  are equivalent 
as quotient maps if there are invertible c.a. $\phi, \psi$ such
that $f=\psi g \phi$. 
We prove Theorem \ref{bignu} below to show that for a c.a. $f$,  
 no property 
defined on equivalence classes of  quotient maps 
can prevent $\nu (f,S_N)$ from being 
arbitrarily close to $N$. To avoid a lengthy digression to 
background, we give a proof assuming familiarity with 
symbolic dynamics; however the statement of Theorem \ref{bignu} 
is self-contained. 
Below,  $h(T)$ denotes the topological entropy of 
$T$.

\begin{lemma}\label{compose} 
Suppose $(\Sigma_A ,S)$ is a mixing shift of finite type (SFT) 
 of positive entropy,  
and $f:\Sigma_A\to \Sigma_A$ is a surjective block code, 
and $\delta >0 $. Then there is an automorphism $\phi$ 
of  $(\Sigma_A,S)$ and a mixing SFT  $(\Sigma ,T)$ contained 
in $S$ such that $h(T)> h(S) - \delta $ and 
the fixed point set of $  \phi f$ contains $\Sigma$. 
\end{lemma}

\begin{proof}
%[Outline of proof]
First, pick a periodic point $y$ in $\Sigma_A$ such that 
$y$ and $f(y)$ have the same least period. (Such $y$ must 
exist: otherwise, $f$ would map the periodic points of prime 
least period to fixed points, and this would imply that 
$f(\Sigma_A)$ is a single point, contradicting surjectivity 
of $f$ and 
 positive entropy of $(\Sigma_A, S)$.) 

Next, e.g. using \cite[Lemma 26.17]{DGS},
find a mixing SFT $(\Sigma_1',T_1')$ in $(\Sigma_A,S)$ 
such that $h(T_1') > h(S) - \delta$ and $y\notin \Sigma_1'$. 
Let $\Sigma_1= f^{-1}(\Sigma_1')$. 
Now easily construct a mixing SFT $(\Sigma_2,T_2)$ with 
$(\Sigma_A,S) \supset (\Sigma_2,T_2)
\supset \Sigma_1 \cup \{y\} $ and 
$\Sigma_2 \cap f^{-1}\{f(y)\}=\{y\}$. 

The restriction of $f$ to $\Sigma_2$ is finite to one 
with degree 1. 
Let $W$ be a magic word for this restriction. 
Let $X_M$ denote the set of points $x$ of
$\Sigma_2$ such that for every $i$ the word $x[i,i+M]$ contains an
occurrence of $W$. Then the restriction of $f$ to $X_M$ will be 
one-to-one, and for large enough $M$ the restriction of $S$ to $X_M$ 
will still have entropy greater than $h(S) -\delta$. 

Now pick $K$ such that for every $n\geq K$, 
$S$ has at least two orbits of length 
$n$ which are not in $T_1$.  Then choose $(\Sigma ,T)$ to be 
an SFT  inside $(X_M,S|X_M)$  such that 
$T$ has no orbits of length smaller than $K$ and also 
such that still $h(T)> h(S) -\delta $. The image of $(\Sigma,T)$ 
under $f$ is an SFT $(\Sigma',T')$ isomorphic to 
$(\Sigma ,T)$, and there is a block 
code $g: \Sigma'\to \Sigma$ such that on $\Sigma$, $gf$ is the identity map.
 By \cite[Theorem 1.5]{BKr2}, 
there is an automorphism $\phi$ of $S$ whose restriction to 
$T'$ equals $g$. Clearly the restriction of 
$ \phi f$ to $T$ is the identity map. 
\end{proof}

\begin{theorem}\label{bignu} 
Suppose $f$ is a surjective c.a. on $N$ symbols 
and $\epsilon >0 $. Then there is an 
invertible c.a. $\phi$ 
 such that $\nu (\phi f,S_N) >  N  - \epsilon$.
\end{theorem}

\begin{proof} 
If $T$ is a mixing shift of finite type with 
$h(T)=\log \lambda$, then 
$\lim_k  |\textnormal{Fix}(T^k)|^{1/k} = \lambda$. 
If this $T$ is a set of fixed points for a c.a. $\psi $ 
on $N$ symbols, it follows that 
$ \nu(\psi ,S_N)\geq \lambda $. 
Now the theorem follows from Lemma \ref{compose}. 
\end{proof}

\begin{remark} 
The statements of Lemma \ref{compose} and Proposition \ref{bignu} 
remain true if $ \phi f$ is replaced by $f \phi $. One 
way to see this is to notice that the systems $(f\phi ,S)$ 
and $(\phi  (f  \phi) \phi^{-1}, \phi S \phi^{-1})=
(\phi  f , S )$ are topologically conjugate. 
\end{remark} 

(2) Now we turn to algebra.  $\Sigma_N$ is 
 a group under coordinatewise addition (mod $N$), and 
some c.a. are group endomorphisms of this group;  
these  are the {\it linear} cellular automata   
whose jointly periodic points were studied in 
\cite{MOW} and later in a number of papers (see 
\cite{CCG,CDN,Jen,LiTh,MST,Su} and their references). 
The algebraic structure allowed a number theoretic description 
of the  way that $f$-periods of jointly periodic points of 
$S_N$ period $n$ vary (irregularly) 
with $n$.  
We show now that when  $f$ is  a linear c.a., 
it is easy to see that $\nu(f,S_N)=\log N$.

\begin{proposition} 
Suppose a c.a. map $f$ is a linear map on $S_N$. 
Then for all large primes $p$, 
$\nu_p(f,S_N)\geq N^{p-1}$. Therefore $\nu(f,S_N)= N$. 
\end{proposition}

\begin{proof}
We use an argument from  the 
proof  of a related result, Proposition 3.2 of 
 \cite{BKi}. 
Let $M$ be the cardinality of the kernel of $f$. 
Suppose $p>M$ and $p$ is prime;   
then $f$ must map orbits of length $p$ 
to orbits of length $p$ (otherwise, some orbit of length 
$p$ would be collapsed to an orbit of length dividing $p$,
i.e. to a fixed point, which would contradict the fact that 
 every point has $M$ preimages, with $M<p$). 
 The fixed points of 
$(S_N)^p$ form a subgroup $H$ which is mapped into itself by $f$, 
and for all $k>0$, the kernel of $f^k$ contains 
no point in an orbit of length $p$, 
so $H\cap \text{ker}(f^k) \subset \text{Fix}(S_N)$. 
For some $k>0$, the restriction of $f^k$ to 
$f^kH$ is injective, and all points in the set 
$f^kH$ are $f$-periodic. Because the kernel of $(f^k)|H$ 
contains at most $N$ points, it follows that at least 
$1/N$ of the fixed points of $S^p$ are periodic for $f$. 
\end{proof}

(3) Algebra can be used in another way.
Frank Rhodes \cite{Rh}, using properties of certain 
families of polynomials presenting c.a. maps, 
 exhibited  a family of noninvertible c.a. $f$ 
 for which there exists $k\in \mathbb N$ such that 
$f$ is injective on $P_{kn}(S_N)$ for all $n\in \mathbb N$. 
Clearly in this case $\nu (f,S_N)=N$. We will not review 
that argument.

(4) We now turn to equicontinuity. A point $x$ is equicontinuous 
for $f$ if for every positive integer $M$ there exists a 
positive integer $K$ such that 
for all points $y$, if $x[-K,K]=y[-K,K]$ then 
$(f^nx)[-M,M]=(f^ny)[-M,M]$ for all $n>0$ . If the surjective c.a. $f$ has $x$
as a point of equicontinuity, and $M$ is chosen larger than the 
span of the block code $f$, and $W$ is the corresponding word 
 $x[-K,K]$, then the following holds: if $z$ is a point 
in which $W$ occurs with bounded gaps, then $z$ is $f$-periodic. 
Thus $\lim_n (1/n)\log \nu_n(f,S_N) = \log (N)$, and moreover 
the convergence is exponentially fast \cite{BT}. 
 Points of equicontinuity may occur in natural examples 
\cite{BlMa, kurka}. 

For many (probably ``most'') surjective 
c.a., the criteria above 
are not applicable. 
This leads to the experimental investigations 
discussed in the next section, and to the possibility raised 
in Questions \ref{Q1} and \ref{Q2} of a general plenitude of jointly 
periodic 
points. Question \ref{Q1} arises 
because in the experimental data, 
the restrictions of the c.a. $f$ to $P_k(S_N)$ are somewhat
reminiscent of a random map on a finite set. 
Since $f$ is a surjective one dimensional c.a. map, there is an $M$ such that 
no point has more than $M$ preimages under $f$. Suppose for 
example $k$ is a prime greater than $M$ and let $\mathcal O_k(S_N)$
denote the set of $S_N$ orbits of size $k$.  Then $f$ defines an  
at most $M$-to-$1$ map $f_k$ from $\mathcal O_k(S_N)$ into itself, 
and we see a possible heuristic: (1) in the absence of 
some additional structure,   
 the sequence $(f_k)$ will reflect  some 
properties of random maps,  and (2) 
an ``additional structure'' such as existence of 
equicontinuity points for $f$  will tend to  produce more rather than 
fewer periodic points. 
The  beautiful and  extensive theory of random maps on finite sets
contains precise asymptotic distributions answering various 
natural questions \cite{Sa}.   
 Here we simply  note that 
for a random map 
on a set of $K$ elements, asymptotically 
 on the order of $\sqrt K$ 
of the elements will lie in cycles 
(whether the map is bounded-to-one  \cite[Theorem 2]{Gr} 
or not  \cite{Sa}), and there will be few big cycles.    

The maps $f_k$ derived from the surjective c.a. $f$
are nonrandom not only in being bounded-to-one, but also 
in that 
 most points have 
the minimal possible number of preimages \cite{hedlund,Ki,LM}. 
To the extent it matters, this seems 
to work in favor of the random maps heuristic behind Question
\ref{Q1}. In particular, it seems that the qualification to the 
random maps analogy offered in \cite[p.252]{MOW}, regarding 
large in-degrees for cellular automata, does not hold for the 
class of surjective c.a.

\section{The maps} \label{maps} 

We examine with our programs 
several cellular automata on $N$ symbols, 
 having or not having various
properties as indicated below. Except for Tables 
\ref{3shiftbiclosdeg9} and \ref{3shift}, all c.a. examined are 
on $N=2$ symbols.  

The c.a. $A$ is the addition map $x_0+x_1$ (mod $N$).  
This c.a. is linear, bipermutative, and everywhere $N$ to one. 
 
The c.a. $B$ is $x_0 +x_1x_2$. This c.a. is left permutative, degree 1,  
not right closing. 

The c.a.  $C$ is $ B\circ B_{rev}$, where 
$B_{rev}=x_0x_1+x_2$.  
This c.a. is degree 1, and it is 
nonclosing, as it is the composition of a not-left-closing
c.a. and a not-right-closing c.a. 

The c.a. $D$ is  the
map $C$ composed with $(S_2)^{-2}$, i.e., 
$D$ is the composition of $x_0 +x_1x_2$ with   $x_{-2}x{-1}+x_0$.
All periodic points for the  golden mean shift (the sequences $x$ in which the 
word $11$ does not occur) become fixed points for $D$ (vs. 
being periodic of varying periods for $C$).

The c.a. $E$ is the composition $A$ followed by $B$. 
This c.a. on $N=2$ symbols  has degree $2$, and is left permutative
but not right closing.

The c.a. $J$ on 2 symbols is $A$ precomposed with the 
automorphism $U$ of $S_2$ which applies  the flip to the symbol 
in the $*$ space of the  frame $1 0 * 1 1$.   
This $U$ is $ x_0 + x_{-2}(1+x_{-1})x_1x_2$, which equals 
$ x_0 + x_{-2}x_1x_2 +  x_{-2}x_{-1}x_1x_2$. 
The c.a. $J$ has degree $N$ and is biclosing, but 
is neither left permutative 
nor right permutative. 

The c.a. $G$ is $ x_{-1}+x_0x_1+x_2$. 
This c.a. on 2 symbols is bipermutative, degree $2$, and is  not linear.  

The c.a. $H$ is the composition $A\circ A \circ U$. 
It has the properties of $F$, except that the degree is now $2^2=4$. 

The c.a. $K$ is the composition $B\circ U$. This c.a. 
is left closing degree 1; it is not left permutative and it is not 
right closing. 

 In addition we use 
a library of surjective 
span 4 and span 5 c.a. due to 
Hedlund, Appel and Welch, 
 who  conducted the early investigation 
\cite{HAW} in which 
they found all surjective c.a. on two symbols of span at most five. 
(This was not trivial, especially in 1963, 
 because there are $2^{32}$ c.a. on two symbols  of span at 
most five.) Among these onto maps of span  four, there are exactly 
32 which are not linear in an end variable (i.e., neither left 
nor right permutative) and which 
send the point $\dots 0000\dots $ to itself. These 32 are 
listed in Table \ref{span4maps}. Any other span four onto 
map which is not linear in an end variable 
is one of these 32 maps $g$ precomposed or postcomposed 
with the flip map $F= x_0+1$. Because 
$gF=F(Fg)F=F^{-1}(Fg)F$, the jointly periodic data for $Fg$ and 
$gF$ will be the same. Altogether, then, we can handle all
surjective span 4 maps not linear in an end variable by examining 
64 maps.

According to \cite{HAW}, there are 141,792 surjective c.a. 
of span 5. These are arranged in \cite{HAW} into classes -- 
linear in end variables, compositions of lower-span maps, and the remainder. 
The remainder class (11,388  maps) is broken down into subclasses 
by patterns of generation, and a less regular 
residual class of 200 maps. These 200 are generated by 26 maps 
\cite[Table XII]{HAW} and various operations. 
We list the codes for this irregular class of 26 maps in 
Table \ref{span5maps}, and use it as a modest 
sample of span 5 maps.

\section{FDense}\label{secfdense} 

The program FDense takes as its input 
a c.a. $f$, an integer $N\geq 2$, a positive integer $m$
and a finite set $\mathcal K$ of positive integers $k$. 
(FDense can also handle sets of maps as inputs, producing 
output for all the maps, and suppressing various data.)  The input $f$ can be given 
by a polynomial or a tabular rule. 
For a given $f$ and  
each $k$ in $\mathcal K$, FDense determines whether the set 
$\text{Per}(f)\cap P_k(S_N)$ is $m$-dense 
(in which case we say that $f$ is $m$-dense at $k$). 
%(i.e., whether 
%every $S_N$-word of length $m$ occurs in some 
%jointly periodic point in  $P_k(S_N)$). 
If not, then FDense will separately list all the $S_N$ 
words of length $m$ which do not appear in any 
periodic point of $f$ in $P_k(S_N)$, in a lexicographically 
truncated form potentially useful for seeing patterns. (For example, 
if $m$ is ten and the word 011 does not occur in the examined points, then 
FDense would list 011* as excluded rather than listing 
all words of length ten beginning with 011.) 

The underlying algorithm for FDense lists all   
words of length $m$ and $k$ in tagged form and operates 
on tags as it moves through the words of length $m$ with 
$f$. 
Memory is the fundamental constraint on FDense. With 
$m$ considerably smaller than $k$, the essential demand on 
memory is the tagged list of $N^k$ words of length $k$. With $N=2$, 
roughly 
$m=13$ and $k=27$ was a practical limit for our machine, 
and this was also quite slow. 
We restricted our investigations almost entirely 
 to the case of $N=2$ symbols for two reasons: 
with $N=2$ we can 
examine longer periods;   
and we would be astonished to find any relation between 
the questions at hand  and $N$. 

The following proposition follows from the data of  
Tables \ref{able13densespan4} and 
\ref{able13densespan4F}.

\begin{proposition}\label{13dense} For every span 4 surjective 
cellular automaton on two symbols, the set of jointly 
periodic points is (at least) 13-dense. 
\end{proposition} 

In Tables \ref{abledensespan4}-\ref{abledenseperm+},   
we applied FDense, for $N=2$ symbols, to check  
 for which $k\leq 24$ 
various other surjective c.a. $f$ are 10-dense at $k$. 

{\it Table \ref{abledensespan4}.} After postcomposition with the 
map $A=x_0+x_1$, the 32 onto span 4 c.a. of Table \ref{span4maps} remain 
10-dense at some $k\leq 24$.  

{\it Table \ref{abledensespan5}.} The 26 irregular span 5 
maps of Table \ref{span5maps} are 10-dense at some $k\leq 24$.  

{\it Table  \ref{abledenseperm+}.} For each of the 32 span 4 maps 
$j$ of  Table \ref{span4maps}, 
 let $p_j(x_0,x_1,x_2,x_3)$ 
denote its defining polynomial.  
Construct a c.a. $f_j$ with defining polynomial 
$x_0+p_j(x_1, 
x_2,x_3,x_4)$. These $f_j$ are demonstrated to be $10$-dense 
at some $k\leq 24$.  

For the c.a. in Tables \ref{abledensespan4}-\ref{abledenseperm+},   
often the least $k$ at which 10-density is achieved lies 
in the range $19-24$. (This is the point of 
Table \ref{abledenseperm+},  as we know already  
from \cite{BKi} that the jointly periodic points of permutative 
c.a. are dense.) 
 This is  consistent 
with the heuristic that apart from possible extra structure the 
c.a. map on points of least period $k$ looks something like 
a random map. For a random map $f$ from a set of $2^k$ points 
into itself, on the order of $\sqrt{2^k}$ points are expected 
to lie in $f$-cycles. For $k=20$, we have $\sqrt{2^k}=2^{10}$. 
(Of course, $10< 24/2$. 
A point of $S_N$-period $20$ will contain up to $20$ distinct 
words of length 10; the words aren't expected to occur  
with complete uniformity;  specific codes are not random.   
For the heuristic of randomness, it is perhaps striking to 
find the rough agreement we do see.) 

 We also checked 10-denseness 
for several c.a. on $2$ symbols with specified properties,  
  described 
in  Section \ref{maps}.

\begin{example}\label{denseplus} [Linear] 
The c.a.  $A = x_0 +x_1$ is 10-dense 
at $k=11,13-24$ out of [10,24].    
\end{example} 

\begin{example}\label{densepermut} [Permutative, not biclosing]
The  c.a. $B$ is 10-dense at 
$k= 22-24$ out of [10,24].  It is 13-dense at only $k=25$ 
out of [13,25]. 
\end{example} 
 
\begin{example}\label{densenonclosing} [Not closing] 
The c.a. $C$ (and likewise $D$)
 is 10-dense at $k=17-24$ 
 out of [1,24], 
and 13-dense for $k=23,24$ out of [13,24]. 
\end{example} 

\begin{example}\label{densebicl} [Degree 2, biclosing, not
permutative]
The c.a. $J$  is  10-dense at 
$k=23-25$ out of [10,25].  It is 13-dense at only $k=25$  
out of [13,25]. 
\end{example} 

In summary, there is reasonable supporting evidence for the Conjecture 
\ref{D}, and the counts seen seem consistent with the random maps 
heuristic.

\section{FPeriod} \label{FPeriod} 

Recall $P_k(S_N)$ denotes the 
set of points fixed by the $k$th power of the full 
shift on $N$ symbols. Each such point $x$ is 
determined by the word $x_0x_1\dots x_{k-1}$. 

The FPeriod program takes as input a c.a. 
 $f$, 
an integer $N\geq 2$ and a finite set of positive 
integers $k$. For each $k$, the program then determines 
data including the following (included in tables cited 
below):  
\begin{itemize} 
\item
$P$:= the number of points in $P_k(S_N)$ which 
are  periodic for $f$. 
\item 
$L$:= the length of the longest $f$-cycle in 
$P_k(S_N)$.  
\end{itemize} 
The program does much more; for the points in 
$P_k(S_N)$, it can produce a complete 
list of $f$ cycle lengths and preperiods 
with multiplicities, and related data such as 
$\nu_k$ and averages. It can also 
do this for points in 
$P^o_k(S_N)$ rather than 
$P_k(S_N)$ (i.e. for points of least shift 
period $k$). The program also has an option 
for producing truncated and 
assembled data for a collection of maps.

The basic algorithm idea of  FPeriod 
is the following. FPeriod takes the given c.a.  $f$ and a given 
shift-period length $k$; stores all $2^k$ words of length $2^k$; and 
then changes various tags on these words as $f$ moves 
through the corresponding periodic points. The tags 
in particular are changed to keep track of how long $f$ 
iterates before returning.  When the program returns to 
a previously visited point, it can deduce the corresponding 
$f$ period and preperiod. 
The essential limit of FPeriod is that for large $k$ 
it becomes 
a horrendous memory hog.  We could conveniently reach period 
$k=23$, and with  patience we could reach $k=25$ or $26$, 
before our memory resources were exhausted. In practice, 
running the program using $N = 2$ and $k = 26$  
 required 1.8 gigabytes  of memory.  

In this section we apply FPeriod to various maps from 
Section \ref{maps} with specific properties, and 
also to many maps of span 4 and 5. The main 
message is that for nonlinear 
maps, we generally see  $\nu_k(f,S_N)$ 
compatible with  affirmative answers to 
Questions \ref{Q1} and \ref{Q2}, and frequently 
the data suggestion strongly that the limit 
$\nu(f,S_N)$ is smaller than $N$.  
Below, unless otherwise indicated, $f$ is defined on 
the full shift $S_N$ with $N=2$, and the symbol set is 
$\{0,1\}$.

{\it Table \ref{ableadditionmap}} [Linear]. 
We exhibit results for the c.a. $A = x_0 +x_1$; here $\nu_k(A,S_2)$ 
is large, consistent with the fact 
$\nu(A,S_2)=2$.

{\it Table \ref{ableF}}  [Biclosing].
We exhibit results for the c.a. $J$, which is $A$ composed 
with an invertible c.a. The composition significantly reduces the 
numbers $\nu_k$.

{\it Table \ref{ableE}} [Linear composed with degree 1 permutative].
We exhibit results for the c.a. $E$. 

{\it Table \ref{ablebipermut}} [Bipermutative].
We exhibit results for the c.a. $G$. 

{\it Table \ref{ableleftpermmap}}  [Permutative, not biclosing].
We exhibit results for the c.a. $B$. 
%degree 1 left permutative map 
%which is not right closing. 

{\it Table \ref{ableK}}   [Closing, not permutative, not biclosing].
We exhibit results for the c.a. $K$. 

{\it  Table \ref{ablenotclosing}}   [Not closing].
We exhibit results for the c.a. $C$. 

{\it Tables \ref{3shiftbiclosdeg9} and 
\ref{3shift}}. We give our only examples 
for a c.a. on more than 2 symbols (they are 
c.a. on 3 symbols).  The pattern is the same 
but we are able to investigate only up to shift period 13. 

{\it Tables \ref{ablespan5a} and \ref{ablespan5b}} [Span 5
irregular]. 
We display data for the 26 irregular maps of span 5 given 
in Table  \ref{span5maps} and discussed in Section \ref{maps}. 

{\it Tables \ref{ablespan4a} and \ref{ablespan4b}} [Span 4]. 
We exhibit data for 
the 32 maps $g$ of Table \ref{span4maps}. 
(This addresses all span 4 surjective c.a. on 2 symbols 
not linear in an end variable, as discussed in Section \ref{maps}.) 

{\it Tables  \ref{ablespan4c} and \ref{ablespan4d}}. [Span 4 composed
with flip].
We exhibit data for the 
 32 maps of Table \ref{span4maps}  postcomposed with the flip
involution
$F=x_0+1$.

{\it Table \ref{ableresolv}} [Permutative comparison].
  $\nu^o_k$ is computed
for  16 left permutative span 5 maps, 
to make a rough comparison of a sample of  maps which are and are not 
linear in an end variable. We see no particular difference.

{\it Table \ref{fppBsamplecycles}.} 
For the map $B$, periods with multiplicity are probed for $k\leq 30$ 
for two samples, of size 10 and size 30. The maximum period is the 
same except for two values of $k$. 

{\it Table \ref{fpBcycles}.} For $B=x_0+x_1x_2$, 
complete data for $B$-periods with 
multiplicity are found by FPeriod (not FProbPeriod) for points 
in $P_k(S_2)$ for $k\leq 22$.

\section{FProbPeriod} \label{secFProbPeriod}

The $k$ for which the 
program FPeriod can explore $f$-periodicity of points in 
$P_k(S_N)$ is limited on account of the memory demands of 
FPeriod. 
This begs for a probablistic approach. For large $k$ 
it  is generally useless to sample points of shift 
period $k$ for $f$-periodicity (commonly, this will be 
a fraction of the shift periodic points exponentially 
small in $k$). Instead, FProbPeriod randomly samples points 
of period $k$ and computes for them the length of the $f$-cycle 
into which they eventually fall. This extends the range of 
$k$ which can be investigated, depending on the map; 
for different maps we've seen practical limits 
at $k=33$ to $k=37$  (typical), to past $50$ (for the linear 
$x_0+x_1$ on two symbols).  In any case, we can search larger 
$k$ than are accessible to us with FPeriod. The program 
FProbPeriod again 
works by listing and tagging, but now only needs to 
keep in memory words for the points visited along an iteration. 
As long as the preperiod and period of the forward orbit 
aren't too large, the program won't crash. 

The input data for FProbPeriod then are the c.a. $f$ ; 
a finite set of periods $k$; the number $N$ of symbols; and 
the number $m$ of points to be randomly sampled for each $k$. 
The program will for each $k$  take  $m$ random 
samples of points from $P_k(S_N )$, and 
find the corresponding periods and preperiods with multiplicity. 
Given $k$, $L$ denotes the largest $f$-period found in the sample. 
For any sequence of samples, 
clearly  $\limsup_k L^{1/k}\leq \limsup_k \nu_k(f,S_N)\leq N$, 
and inequalities must become sharp in some cases 
($f$ linear or $f$ of finite order).  
Still, the data we see seems consistent 
with positive answers to Questions \ref{Q1} and \ref{Q2}.  

%$\lim_k L^{1/k}=\sqrt N$ and $\lim_k \nu_k(f,S_N)= \nu(f,S_N)=N$. 
%For example, if  $f^M=\text{Id}$ for some $M>0$, then $f$-orbit 
%lengths are bounded and all points are $f$-periodic. 
%Also, recall that for a linear surjective c.a., the 
%maximum $f$-period seen in $P_k(S_N)$ has $\limsup_k$ equal to 
%$\sqrt N$. 

The specific maps cited below are described in Section \ref{maps}. 

{\it Table \ref{fppBcycles}.}  
For sample size $m=10$, for the (degree one, left permutative,
not right closing) map $B=x_0 +x_1x_2$, the (eventual) periods are
listed with their multiplicities  in the sample, 
for $1\leq k \leq 37$.  

{\it Table \ref{fppBsamplecycles}.} 
For the map $B$, periods with multiplicity are probed for $k\leq 30$ 
for two samples, of size 10 and size 30. The maximum period is the 
same except for two values of $k$. By comparison with 
Table \ref{fpBcycles}, one sees that the 
size 30 sample in Table
\ref{fppBsamplecycles} found the largest period except at $k=12$ 
(where it found period 56 but not the maximum period 60).

{\it Table \ref{fppAsamplecycles}.} 
For the linear c.a.  $A$, periods with multiplicity are 
probed for $k\leq 49$ for two sample sizes, 10 and 30. 
The results are almost identical. 
 
{\it Table \ref{fpptable}.} 
For sample size $m=10$, for $1\leq k \leq 37$, 
the numbers $L^{1/k}$ are computed for several c.a. described in 
Section \ref{maps}:
$A,B,C,E,G,H,J,K$. The corresponding preperiod 
data is displayed in {\it Table \ref{fppptable}.} 

{\it Table \ref{fppCcycles}.} 
For sample size $m=10$, for $1\leq k \leq 32$, 
 the sampled periods for the nonclosing c.a. $C$ are 
listed with their multiplicities  in the sample.  

{\it Table \ref{fppDcycles}.} 
For sample size $m=10$, for $1\leq k \leq 32$, 
 the sampled periods for the nonclosing c.a.  $D$ are 
listed with their multiplicities  in the sample.  

{\it Table \ref{preperB}.} 
This table lists the preperiods found for $B$ by FProbPeriod 
for the sample size 10 in the range $18\leq k\leq 35$. 

{\it Table \ref{preperC}.} 
This table lists the preperiods found for $C$ by FProbPeriod 
for the sample size 10 in the range $18\leq k\leq 35$.

%\section{Discussion of the programs} \label{program}
%just a bit in the intro. 

\clearpage		   
\appendix

\section{Tables of some span 4 and 5 c.a.} 

\begin{table}[ht]
\begin{center}
\begin{tabular}{|c|c|c|c|}
\hline
Map& Tabular rule         &  Map& Tabular rule        \\
\hline 
1  &  0000 1111 0010 1101 & 17 &  0011 1001 1100 1100  \\
2  &  0000 1111 0100 1011 & 18 &  0011 1010 0011 1100  \\
\hline 
3  &  0001 1100 0011 1110 & 19 &  0011 1010 1100 0011  \\
4  &  0001 1110 0101 1010 & 20 &  0011 1100 0101 0011  \\
\hline 
5  &  0010 1001 0110 1101 & 21 &  0011 1100 0101 1100  \\
6  &  0010 1101 0000 1111 & 22 &  0011 1100 1010 0011  \\
\hline 
7  &  0011 0011 0110 0011 & 23 &  0011 1100 1010 1100  \\
8  &  0011 0011 0110 1100 & 24 &  0011 1110 0001 1100  \\
\hline 
9  &  0011 0011 1001 0011 & 25 &  0100 1001 0110 1011  \\
10 &  0011 0011 1001 1100 & 26 &  0100 1011 0000 1111  \\
\hline 
11 &  0011 0101 0011 1100 & 27 &  0101 1010 0001 1110  \\
12 &  0011 0101 1100 0011 & 28 &  0101 1010 0111 1000  \\
\hline
13 &  0011 0110 0011 0011 & 29 &  0110 1011 0100 1001  \\
14 &  0011 0110 1100 1100 & 30 &  0110 1101 0010 1001  \\
\hline 
15 &  0011 1000 0111 1100 & 31 &  0111 1000 0101 1010  \\
16 &  0011 1001 0011 0011 & 32 &  0111 1100 0011 1000  \\
\hline
\end{tabular}
\vskip .1in 
\caption{The 32 span 4 onto c.a. of the 2 shift 
which fix $\dots 000\dots$ and are 
not linear in an end variable 
\cite[Table I]{HAW}. Maps 2, 6, 7 and 16 are 
one-to-one. The rule above for map 30  corrects a misprint in 
\cite[Table I]{HAW}.}\label{span4maps}
\end{center}
\end{table}

\begin{table}[ht]
\begin{center}
\begin{tabular}{|c|c|c|c|}
\hline
Map& Tabular rule         &  Map& Tabular rule        \\
\hline
 1 &  0001 0111 1110   1000  0001  0111 1111 0000& 14 &  0100 1101 1111   0000 0100  1101 1011 0010 \\ 
 2 &  0001 1011 0111   0100  1110  0100 1111 0000& 15 &  0110 0001 1010   1011 0110  0001 0110 0111 \\
\hline			     
 3 &  0010 0010 1111   0011  0010  1110 0000 1111& 16 &  0110 1000 0111   1001 0110  0001 1110 1001 \\
 4 &  0010 1001 0110   1101  0100  1001 0110 1011& 17 &  0110 1011 1100   0010 0100  1011 0001 1101 \\
\hline			     
 5 &  0010 1110 0000   1111  0010  1110 1111 0000& 18 &  0111 0001 1011   0010 0111  0001 1000 1110 \\
 6 &  0100 0111 0001   0111  1011  1000 0000 1111& 19 &  0111 0010 1011   0100 0111  0010 0111 1000 \\
\hline			     
 7 &  0100 0111 0100   1011  1000  1011 0100 1011& 20 &  0111 1000 0100   1011 0111  1000 0111 1000 \\
 8 &  0100 1011 1000   0111  0100  1011 0100 1011& 21 &  0111 1000 0100   1011 0111  1000 1011 0100 \\ 
\hline			     
 9 &  0100 1101 1011   0010  1000  1110 1011 0010& 22 &  0111 1000 0100   1011 0111  1000 1111 0000 \\
10 &  0100 1101 1011   0010  1100  1100 1011 0010& 23 &  0111 1000 0100   1101 0111  1000 1000 1110 \\
\hline			     
11 &  0100 1101 1101   0010  0011  0011 1101 0010& 24 &  0111 1011 1000   0100 0100  1011 0000 1111 \\
12 &  0100 1101 1101   0010  0111  0001 1101 0010& 25 &  0111 1011 1100   0000 0100  1011 0000 1111 \\
13 &  0100 1101 1101   0010  1111  0000 1101 0010& 26 &  0111 1011 1100   0000 0100  1011 0100 1011 \\
\hline
\end{tabular}
\vskip .1in
\caption{\label{span5maps} 26 irregular span 5 onto maps 
of the 2 shift 
which fix $\dots 000\dots$ and are 
not linear in an end variable 
\cite[Table XII]{HAW}.}
\end{center}
\end{table}

\newpage 

\section{FDense Tables} \label{fdensetables}

\begin{table}[ht]                                 
\begin{center}                                
\begin{tabular}{|c|c|c||c|c|c|c|c|}           
\hline                                        
Map& 10-dense at   & 13-dense at & Map& 10-dense at  & 13-dense at \\    
\hline                                                                   
1  &    11-24      &    13-24    & 17 & 17,18,20-24  & 24     \\         
2  &    10-24      &    13-24    & 18 & 17,19-24     & 23-24  \\         
\hline                                       
3  &    18-24      &    24       & 19 & 19-24        & 23-24  \\ 
4  &    21-24      &    (27)     & 20 & 17,19-23     & (25)   \\ 
\hline                                                                           	
5  &    17,19-23   &    (25     & 21 & 19-24        & 23-24  \\                                
6  &    10-24      &    13-24    & 22 & 19-24        & 23-24   \\                             
\hline                                                	
7  &    10-24      &    13-24    & 23 & 19-24        & 21-22  \\               
8  &    21-24      &    (27)     & 24 & 17,19-24     & 23-24  \\               
\hline                                                	
9  &    11-24      &    13-24    & 25 & 17,19-23     & (25)    \\       
10 &    19,21-24   &    24       & 26 & 11-24        & 13-24   \\            
\hline                                                	
11 &    18-24      &    24       & 27 & 19, 21-24    & 24     \\     
12 &    17,19-23   &    (25)     & 28 & 22-24        & (25)   \\        
\hline                                                	
13 &    11-24      &    13-24    & 29 & 19-24        & 23-24  \\              
14 &    22-24      &    (25)     & 30 & 19-24        & 23-24  \\               
\hline                                                	
15 &    19-24      &    23-24    & 31 & 17,18,20-24  & 24     \\          
16 &    10-24      &    13-24    & 32 & 19-24        & 21-22   \\      
\hline		     	      
\end{tabular}	     	
\vskip .1in
\caption{The map numbers refer to the 32 span 4 maps of Table
\ref{span4maps}. Table \ref{able13densespan4} shows for the given sample
of maps, and for $m=10$ and $m=13$, for which $k$ in the range 
$[m,24]$ 
 the jointly periodic points 
in $P_k(S_2)$ are $m$-dense. If the map is not $m$-dense
for any $k$ in this range, then the number listed in parentheses is 
the smallest $k$ for which the jointly periodic points 
in $P_k(S_2)$ are $m$-dense. 
}\label{able13densespan4}
\end{center}   
\end{table}

\begin{table}[ht]                                 
\begin{center}                                
\begin{tabular}{|| c|c|c || c|c|c ||}           
\hline                                        
Map& 10-dense at   & 13-dense at & Map& 10-dense at  & 13-dense at \\    
\hline                                                                   
F$\circ$1  &  11-24      & 13,15-24  	& F$\circ$17 & 19,21-24    &         24   \\
F$\circ$2  &  11-24      & 13-24  	& F$\circ$18 & 19-24       &  23-24       \\
\hline                       				       			  
F$\circ$3  &  19-24      & 21-22  	& F$\circ$19 & 19-24       & (27)   	  \\
F$\circ$4  &  22-24      & (25)  	& F$\circ$20 & 21-24       & 23   	  \\
\hline                       							  			     
F$\circ$5  &  21-24      & 23  	& F$\circ$21 & 17,19-24    &23-24             	  \\
F$\circ$6  &  10-24      & 13-24  	& F$\circ$22 & 19-24       & (27) 	  \\
\hline                       				       			  
F$\circ$7  &  20-24      & 13-24  	& F$\circ$23 & 18-24       &24       	  \\
F$\circ$8  &  22-24      & (25)  	& F$\circ$24 & 19-24       &23-24         \\
\hline                       				       			  
F$\circ$9  &  11-24      &13,15-24   	& F$\circ$25 & 21-24       &23      	  \\
F$\circ$10 &  17-18,21-24&24   	& F$\circ$26 & 11-24       &13,15-24         	  \\
\hline                       				       			  
F$\circ$11 &  19-24      &21,22   	& F$\circ$27 & 17-18,20-24 &24       	  \\
F$\circ$12 &  21-24      &23   	& F$\circ$28 & 21-24       &(27)      		  \\
\hline                       				       			  
F$\circ$13 &  11-24      &13,15-24   	& F$\circ$29 & 19-24       &(27)      	  \\
F$\circ$14 &  21-24      &(27)  	& F$\circ$30 & 19-24       &(27)         \\
\hline                       				       			  
F$\circ$15 &  17-19,20-24&23-24   	& F$\circ$31 & 19,21-24    &24       	  \\
F$\circ$16 &  10-24      & 13-24  	& F$\circ$32 & 18-24 &24                  \\
\hline		     	      
\end{tabular}	     	
\vskip .1in
\caption{The map numbers refer to the 32 span 4 maps of Table
\ref{span4maps}. $F$ is the involution $F=x_0+1$. 
Table \ref{able13densespan4F} shows for the given sample
of maps, and for $m=10$ and $m=13$, for which $k$ in the range 
$[m,24]$ the jointly periodic points 
in $P_k(S_2)$  $m$-dense. If the map is not $m$-dense
for any $k$ in this range, then the number listed in parentheses is 
the smallest $k$ for which the jointly periodic points 
in $P_k(S_2)$ are $m$-dense. 
}\label{able13densespan4F}
\end{center}   
\end{table}

\begin{table}[ht]
\begin{center}
\begin{tabular}{||c|c||c|c||c|c||c|c|} 
\hline
Map& 10-dense at     & Map& 10-dense at  &   Map& 10-dense at     &  Map& 10-dense at   \\      
\hline                                 	 
F$\circ$1  &  11-24      & F$\circ$17 & 19,21-24    &   D$\circ$1  &  20-24     & D$\circ$17 & 21-24   \\        
F$\circ$2  &  11-24      & F$\circ$18 & 19-24       &   D$\circ$2  &  19,22-24  & D$\circ$18 & 21-24     \\          
\hline                                                       		   		       
F$\circ$3  &  19-24      & F$\circ$19 & 19-24       &   D$\circ$3  &  20-24     & D$\circ$19 & 22-24   \\         
F$\circ$4  &  22-24      & F$\circ$20 & 21-24       &   D$\circ$4  &  21,23-24  & D$\circ$20 & 20,24   \\         
\hline                                                		   		                                                           	      
F$\circ$5  &  21-24      & F$\circ$21 & 17,19-24    &   D$\circ$5  &  20,24     & D$\circ$21 & 21,23-24\\                 
F$\circ$6  &  10-24      & F$\circ$22 & 19-24       &   D$\circ$6  &  20,22,24&D$\circ$22  & 24 \\       
\hline                                                       		   		       
F$\circ$7  &  20-24      & F$\circ$23 & 18-24       &   D$\circ$7  &  19,22-24  &D$\circ$23 & 21-24 \\               
F$\circ$8  &  22-24      & F$\circ$24 & 19-24       &   D$\circ$8  &  21,23-24  & D$\circ$24 &21-24     \\               
\hline                                                       		   		       
F$\circ$9  &  11-24      & F$\circ$25 & 21-24       &   D$\circ$9  &  22-24     & D$\circ$25 &22-24        \\       
F$\circ$10 &  17-18,21-24& F$\circ$26 & 11-24       &   D$\circ$10 &  21,23-24  & D$\circ$26 &22-24      \\            
\hline                                                       		   		       
F$\circ$11 &  19-24      & F$\circ$27 & 17-18,20-24 &   D$\circ$11 &  22-24     &D$\circ$27 & 21,23-24        \\     
F$\circ$12 &  21-24      & F$\circ$28 & 21-24       &   D$\circ$12 & 20-21,23-24& D$\circ$28 &22-24       \\        
\hline                                                       		   		       
F$\circ$13 &  11-24      & F$\circ$29 & 19-24       &   D$\circ$13 &  20-24     & D$\circ$29 &22-24  \\              
F$\circ$14 &  21-24      & F$\circ$30 & 19-24       &   D$\circ$14 &  22-24     & D$\circ$30 &24 \\               
\hline                                                       		   		       
F$\circ$15 &  17-19,20-24& F$\circ$31 & 19,21-24    &   D$\circ$15 &  21,23-24  & D$\circ$31 &21-24      \\          
F$\circ$16 &  10-24      & F$\circ$32 & 18-24       &   D$\circ$16 &  20,22,24  & D$\circ$32 &21-24      \\      
\hline		     	      
\end{tabular}	     	
\vskip .1in  
\caption{\label{abledensespan4}
The c.a. listed are compositions, e.g. 
$D\circ j$ is  map $j$ followed 
by $D$.
The map numbers $j$ refer to the 32 span 4 maps of Table
\ref{span4maps}. The map $D$ is given by
$x_0+x_1$. The map $F$ is the 
flip involution $F=1+x_0$.
The data on $F\circ j$ are copied in from Table 
\ref{able13densespan4F} for contrast with $D\circ j$. 
}
\end{center}  
\end{table}

\begin{table}[ht]
\begin{center}
\begin{tabular}{||c|c||c|c||c|c||c|c||}
\hline
j& 10-dense at & j  & 10-dense at   & j  & 10-dense at& j  & 10-dense at \\      
\hline                                                    	  
1  & 18-24        & 9  & 17-19,21-24& 17 & 18-20,22-24 & 25 & 11,13-24   \\        
2  & 19,21-24     & 10 & 19-24      & 18 & 17,21-24    & 26 & 11,13-24   \\          
\hline         	    	      	   
3  & 17-24        & 11 & 17,19-24   & 19 & 17-24       &  &  \\         
4  & 18-19,21-24  & 12 & 18-24      & 20 & 18-24 &  &    \\         
\hline         	    	      	   
5  & 17-24        & 13 & 16-23      & 21 & 19-24 &  &    \\                 
6  & 15-24        & 14 & 16-24      & 22 & 19-24 &  &   \\              
\hline         	    	      	   
7  & 16-24        & 15 & 18-24      & 23 & 20-24 &  &    \\   
8  & 11,13-24     & 16 & 18-19,21-24& 24 & 11,13-24   &  &    \\             
\hline		     	      
\end{tabular}	     	
\vskip .1in
\caption{ 
The map numbers refer to the 26 ``irregular''
span 5  maps of  \cite{HAW}[Table XII], 
copied in    Table \ref{span5maps}. 
\ref{span4maps}. 
}\label{abledensespan5}
\end{center}   
\end{table}

\begin{table}[ht]
\begin{center}
\begin{tabular}{||c|c||c|c||c|c||c|c||}
\hline
j& 10-dense at & j  & 10-dense at   & j  & 10-dense at& j  & 10-dense at \\      
\hline                                                    	  
1  & 20-24     & 9  & 22-23         & 17 & 21,24 & 25 & 20,23-24   \\        
2  & 18-24     & 10 & 21,24         & 18 & 22,24 & 26 & 19,22-24   \\          
\hline           		    		   
3  & 20,22-24  & 11 & 21-24         & 19 & 23,24 & 27 & 21,24    \\         
4  & 21,23-24  & 12 & 21-24         & 20 & 22-24 & 28 & 21,23-24   \\         
\hline           		    		   
5  & 21,23-24  & 13 & 21-24         & 21 & 22,24 & 29 & 20,22-23   \\                 
6  & 19-24     & 14 & 21-24         & 22 & 23,24 & 30 & 19-24     \\              
\hline           		    		   
7  & 20,22-24  & 15 & 21-24         & 23 & 21-24 & 31 & 21,24      \\   
8  & 21-24     & 16 & 20-21,23      & 24 & 21-24 & 32 & 20,22-24   \\             
\hline		     	      
\end{tabular}	     	
\vskip .1in
\caption{As in Table \ref{abledensespan4},
the map numbers refer to the 32 span 4 maps of Table
\ref{span4maps}. 
For such a map $j$, let $p_j(x_0,x_1,x_2,x_3)$ 
be the polynomial such that $(jx)_0= p_j(x_0,x_1,x_2,x_3)$.  
Then a row $j$ of Table \ref{abledenseperm+} 
 refers to the map $f_j$ such that 
 $(f_jx)_0=x_0+p_j(x_1, 
x_2,x_3,x_4)$. Equivalently, $f_j= x_0+ (j\circ S_2)$.  
For the map $f_j$, all $k$ in the range $[10,24]$ at which 
$f_j$ is $10$-dense are listed. 
\newline 
By comparison to Table \ref{abledensespan4}, we see that altering 
the rules $j$ as we have produces maps which look more ``random'' 
in the sense that all the $10$-density is being achieved at $k$
on the order of $2\times 10$. Note that the maps in 
Table \ref{abledenseperm+} are by construction left permutative, 
and for these the jointly periodic points are known to be dense
\cite{BKi}.
}\label{abledenseperm+}
\end{center}   
\end{table}			   
\clearpage

\section{FPeriod Tables}\label{fperiodtables}  

Output obtained from the FPeriod Program, discussed 
in Section \ref{FPeriod}, is compiled in the Tables below. 
In Tables \ref{ableadditionmap}-\ref{3shift},
 for a c.a. $f$,  the numbers given 
for a  row $k$ are computed with respect to $P_k(S_N)$, 
the set of points fixed by the shift $S_N$. 
Except for Tables \ref{3shiftbiclosdeg9} and \ref{3shift}, 
the number of symbols $N$ is $2$.  
For a given $k$, 
$L$ denotes the maximum $f$-period of a point in $P_k(S_N)$;
$P$ denotes the number of points in $P_k(S_N)$ which are $f$-periodic 
(so, $P + \text{Not-}P=N^k$); 
and $\nu_k$ denotes $\nu_k(f,S_N)=P^{1/k}$.   
In some later tables, $\nu^o_k$ is used to 
denote the $k$th root of the number 
of points of {\it least} $S_n$-period $k$ which are periodic for $f$. 
 The preperiod of a 
point $x$ is the smallest nonnegative integer $j$ such that 
$f^j(x)$ is $f$-periodic. 

\begin{table}[ht]
\begin{center}
\begin{tabular}{|r|r|r|r| r|r|r|r| r|r|}
\hline 
    &  Fraction  &  &          &    &  &        & Average   & Average     &  Maximum  \\
$k$ &  Periodic  & $\nu_k$     & $L^{1/k}$ & P  &L & Not-P    & Period  &  Preperiod     & Preperiod \\
\hline 
1   &  0.500000  & 1.00    & 1.00    &	1          & 1        &	  1        &   1.00       &    0.50       &    1 \\    
2   &  0.250000  & 1.00    & 1.00    &	1          & 1        &	  3        &   1.00       &    1.25       &    2 \\        
3   &  0.500000  & 1.58    & 1.44    &	4          & 3        &	  4        &   2.50       &    0.50       &    1 \\        
4   &  0.062500  & 1.00    & 1.00    &	1          & 1        &	  15       &   1.00       &    3.06       &    4 \\        
5   &  0.500000  & 1.74    & 1.71    &	16         & 15       &	  16       &   14.12      &    0.50       &    1 \\    
\hline		       		 										
6   &  0.250000  & 1.58    & 1.34    &	16         & 6        &	  48       &   5.12       &    1.25       &    2 \\    
7   &  0.500000  & 1.81    & 1.32    &	64         & 7        &	  64       &   6.91       &    0.50       &    1 \\    
8   &  0.003906  & 1.00    & 1.00    &	1          & 1        &	  255      &   1.00       &    7.00       &    8 \\    
9   &  0.500000  & 1.85    & 1.58    &	256        & 63&	  256      &   62.05      &    0.50       &    1 \\     
10  &  0.250000  & 1.74    & 1.40    &	256        & 30       &	  768      &   29.01      &    1.25       &    2 \\    
\hline		       		 										
11  &  0.500000  & 1.87    & 1.69    &	1,024       & 341      &	  1024     &   340.67     &    0.50       &    1 \\    
12  &  0.062500  & 1.58    & 1.23    &	256        & 12       &	  3840     &   11.57      &    3.06       &    4 \\    
13  &  0.500000  & 1.89    & 1.67    &	4,096       & 819      &	  4096     &   818.80     &    0.50       &    1 \\    
14  &  0.250000  & 1.81    & 1.20    &	4,096       & 14       &	  12,288    &   13.89      &    1.25       &    2 \\    
15  &  0.500000  & 1.90    & 1.19    &	16,384      & 15       &	  16,384    &   14.99      &    0.50       &    1 \\    
\hline		       		 										
16  &  0.000015  & 1.00    & 1.00    &	1          & 1        &	  65535    &   1.00       &    15.00      &   16\\    
17  &  0.500000  & 1.92    & 1.38    &	65,536      & 255      &	  65,536    &   254.33     &    0.50       &    1 \\    
18  &  0.250000  & 1.85    & 1.30    &	65,536      & 126      &	  196,608   &   125.73     &    1.25       &    2 \\    
19  &  0.500000  & 1.92    & 1.62    &	262,144     & 9,709     &	  262,144   &   9708.96    &    0.50       &    1 \\    
20  &  0.062500  & 1.74    & 1.22    &	65,536      & 60       &	  983,040   &   59.88      &    3.06       &    4 \\    
\hline 		       		 										
21  &  0.500000  & 1.93    & 1.21    &	1,048,576    & 63       &	  1,048,576  &   62.99      &    0.50       &    1 \\    
22  &  0.250000  & 1.87    & 1.34    &	1,048,576    & 682      &  3,145,728  &   681.67     &    1.25       &    2 \\    
23  &  0.500000  & 1.94    & 1.39    &	4,194,304    & 2,047     &	  4,194,304  &   2047.00    &    0.50       &    1 \\    
\hline 
\end{tabular}
\vskip .1in
%this table was abstracted from the output of program inA
\caption{The c.a. is $A=x_0+x_1$ on the 2-shift: a linear, two-to-one
map. If $k=q2^j$ with $q$ odd and $j$ a nonnegative integer, 
then there are exactly $q2^{-(j+1)}$ points in $P_k(S_2)$
which are $f$-periodic. Thus,    $\nu_k=2^{(k-1)/k}$
if $k$ is odd. 
}\label{ableadditionmap}
\end{center}
\end{table}

\begin{table}[ht]
\begin{center}
\begin{tabular}{|r|r|r|r| r|r|r|r| r|r|}
\hline 
    &  Fraction  &   &          &  & &        &     Average   & Average     &  Maximum  \\
$k$ &  Periodic  & $\nu_k$        &  $L^{1/k}$& P         &L & Not-P   &    Period     &  Preperiod     & Preperiod \\ 
\hline 
1   &   .5000  & 1.00   &  1.00   &  1        &   1      &     1        &   1.00      &     0.50    &  1        \\
2   &   .2500  & 1.00   &  1.00   &  1        &   1      &     3        &   1.00      &     1.25    &  2        \\
3   &   .5000  & 1.58   &  1.44   &  4        &   3      &     4        &   2.50      &     0.50    &  1        \\
4   &   .3125  & 1.49   &  1.41   &  5        &   4      &     11       &   2.50      &     1.31    &  3        \\
5   &   .3437  & 1.61   &  1.58   &  11       &   10     &     21    &   9.44      &     0.97    &  2       \\
\hline 
6   &   .4375  & 1.74   &  1.61   &  28       &   18     &     36       &   11.31     &     0.69    &  2       \\
7   &   .0625  & 1.34   &  1.32   &  8        &   7      &     120      &   6.91      &     4.00    &  7        \\
8   &   .0195  & 1.22   &  1.18   &  5        &   4      &     251      &   3.25      &     6.58    &  12       \\
9   &   .1484  & 1.61   &  1.58   &  76       &   63     &     436      &   58.26     &     3.17    &  7        \\
10  &   .0888  & 1.57   &  1.52   &  91       &   70     &     933      &   18.17     &     7.77    &  17       \\
\hline 
11  &   .0703  & 1.57   &  1.46   &  144      &   66     &     1,904     &   65.35     &     5.52    &  14       \\
12  &   .0576  & 1.57   &  1.36   &  236      &   42     &     3,860     &   24.44     &     10.98   &  34      \\
13  &   .0350  & 1.54   &  1.53   &  287      &   273    &     7,905     &   217.65    &     11.93   &  29       \\
14  &   .0201  & 1.51   &  1.39   &  330      &   105    &     16,054    &   12.65     &     36.60   &  74       \\
15  &   .0123  & 1.49   &  1.44   &  404      &   255    &     32,364    &   179.68    &     35.36   &  91       \\
\hline 
16  &   .0232  & 1.58   &  1.54   &  1,525     &   1,008   &     64,011    &   272.23    &     33.28   &  98       \\
17  &   .0286  & 1.62   &  1.52   &  3,758     &   1,377   &     127,314   &   913.23    &     31.04   &  114      \\
18  &   .0091  & 1.54   &  1.53   &  2,386     &   2,250   &     259,758   &   2,026.85   &     55.23   &  152     \\
19  &   .0039  & 1.49   &  1.47   &  2,091     &   1,672   &     522,197   &   1,658.11   &     91.44   &  251      \\
20  &   .0015  & 1.44   &  1.31   &  1,635     &   240    &     1,046,941  &   14.16     &     279.12  &  575      \\
\hline
21  &   .0046  & 1.54   &  1.48   &  9,650     &   4,326   &     2,087,502  &   461.24    &     244.11  &  638      \\
22  &   .0011  & 1.47   &  1.40   &  4,896     &   1,848   &     4,189,408  &   1,158.45   &     274.42  &  647      \\
23  &   .0027  & 1.54   &  1.53   &  23,461    &   19,297  &     8,365,147  &   18,849.71  &     269.70  &  824      \\
\hline
\end{tabular}
\vskip .1in
%this map was from the output of program inF, from file A#U.txt
\caption{The c.a. is $J$, the composition $x_0+x_1$ followed by 
the involution $U=x_0 + x_{-2}x_1x_2 +  x_{-2}x_{-1}x_1x_2$.
This invertible c.a.  $U$ is 
the involution of the 2-shift which replaces 
$x_0$ with $x_0+1$ when 
$x[-2,2]=10x_011$. $J$ is biclosing but not permutative. 
} 
\label{ableF}
\end{center}
\end{table}

\begin{table}[ht]
\begin{center}
\begin{tabular}{|r|r|r|r| r|r|r|r| r|r|}
\hline 
    &  Fraction  &  &          &  & &        &     Average   & Average     &  Maximum  \\
$k$ &  Periodic  & $\nu_k$          &  $L^{1/k}$  & P         &L & Not-P   &    Period     &  Preperiod     & Preperiod \\ 
\hline 
1   & 0.5000  & 1.00   &  1.00    & 1        &   1      &     1       &    1.00      &     0.50     &      1        \\
2   & 0.2500  & 1.00   &  1.00    & 1        &   1      &     3       &    1.00      &     0.75     &      1        \\
3   & 0.1250  & 1.00   &  1.00    & 1        &   1      &     7       &    1.00      &     1.62     &      2        \\
4   & 0.3125  & 1.49   &  1.41    & 5        &   4      &     11      &    2.50      &     0.94     &      2        \\
5   & 0.3437  & 1.61   &  1.58    & 11       &   10     &     21      &    9.44      &     0.66     &      1       \\
\hline 
6   & 0.0156  & 1.00   &  1.00    & 1        &   1      &     63      &    1.00      &     3.52     &      5       \\
7   & 0.2812  & 1.66   &  1.54    & 36       &   21     &     92      &    17.62     &     1.05     &      2        \\
8   & 0.0195  & 1.22   &  1.18    & 5        &   4      &     251     &    3.91      &     5.65     &      11       \\
9   & 0.0546  & 1.44   &  1.22    & 28       &   6      &     484     &    5.82      &     5.18     &      11       \\
10  & 0.1767  & 1.68   &  1.46    & 181      &   45     &     843     &    29.75     &     2.14     &      7        \\
\hline 
11  & 0.0703  & 1.57   &  1.55    & 144      &   132    &     1904    &    98.08     &     6.76     &      19       \\
12  & 0.0012  & 1.14   &  1.12    & 5        &   4      &     4091    &    1.01      &     19.25    &      36      \\
13  & 0.0556  & 1.60   &  1.49    & 456      &   182    &     7736    &    162.94    &     18.54    &      49       \\
14  & 0.0261  & 1.54   &  1.35    & 428      &   70     &     15956   &    28.54     &     18.35    &      55       \\
15  & 0.0342  & 1.59   &  1.45    & 1121     &   285    &     31647   &    138.58    &     21.60    &      58       \\
\hline 
16  & 0.0074  & 1.47   &  1.47    & 485      &   480    &     65051   &    430.96    &     71.09    &      146      \\
17  & 0.0160  & 1.56   &  1.55    & 2109     &   1734   &     128963  &    1633.83   &     51.36    &      169      \\
18  & 0.0060  & 1.50   &  1.41    & 1594     &   549    &     260550  &    334.44    &     70.40    &      233     \\
19  & 0.0046  & 1.50   &  1.45    & 2452     &   1197   &     521836  &    834.45    &     92.00    &      227      \\
20  & 0.0058  & 1.54   &  1.50    & 6165     &   3640   &     1042411 &    2700.37   &     70.21    &      211      \\
\hline 
21  & 0.0017  & 1.47   &  1.36    & 3627     &   693    &     2093525 &    585.86    &     356.39   &      817      \\
22  & 0.0033  & 1.54   &  1.46    & 14004    &   4147   &     4180300 &    3305.59   &     251.62   &      864      \\
23  & 0.0022  & 1.53   &  1.53    & 18746    &   18538  &     8369862 &    18491.96  &     262.30   &      900      \\
\hline
\end{tabular}
\vskip .1in
%this map was from the output of program inE, from file B#A.txt
\caption{The c.a. $E$ is the composition $x_0+x_1$ followed by
$x_0+x_1x_2$, a linear 2-to-1 map followed by a 
degree 1 left permutative map.}
\label{ableE}
\end{center}
\end{table}

\begin{table}[ht]
\begin{center}
\begin{tabular}{|r|r|r|r| r|r|r|r| r|r|}
\hline 
    &  Fraction  &   &           &  & &        &     Average   & Average     &  Maximum  \\
$k$ &  Periodic  & $\nu_k$    &  $L^{1/k}$     & P         &L & Not-P   &    Period     &  Preperiod     & Preperiod \\ 
\hline 
1   &   1.0000 &  2.00  &   1.00   &  2        &   1    &       0       &    1.00    &       0.00   & 0       \\
2   &   .5000 &  1.41  &   1.00   &  2        &   1    &       2       &    1.00    &       0.50   & 1       \\
3   &   .2500 &  1.25  &   1.00   &  2        &   1    &       6       &    1.00    &       1.12   & 2       \\
4   &   .1250 &  1.18  &   1.00   &  2        &   1    &       14      &    1.00    &       1.62   & 3       \\
5   &   .2187 &  1.47  &   1.37   &  7        &   5    &       25      &    1.62    &       2.03   & 4      \\
\hline 
6   &   .4062 &  1.72  &   1.51   &  26       &   12   &       38      &    4.56    &       1.20   & 3      \\
7   &   .0703 &  1.36  &   1.32   &  9        &   7    &       119     &    6.91    &       4.65   & 9       \\
8   &   .0703 &  1.43  &   1.41   &  18       &   16   &       238     &    1.94    &       6.98   & 12      \\
9   &   .1796 &  1.65  &   1.48   &  92       &   36   &       420     &    15.52   &       2.55   & 7       \\
10  &   .0263 &  1.39  &   1.17   &  27       &   5    &       997     &    1.07    &       9.08   & 15      \\
\hline 
11  &   .1782 &  1.70  &   1.53   &  365      &   110  &       1,683    &    77.16   &       3.79   & 16      \\
12  &   .0122 &  1.38  &   1.30   &  50       &   24   &       4,046    &    17.51   &       10.57  & 26     \\
13  &   .1049 &  1.68  &   1.53   &  860      &   260  &       7,332    &    199.90  &       6.20   & 21      \\
14  &   .0056 &  1.38  &   1.37   &  93       &   84   &       16,291   &    70.69   &       22.86  & 48      \\
15  &   .0340 &  1.59  &   1.43   &  1,117     &   225  &       31,651   &    117.64  &       13.52  & 42      \\
\hline 
16  &   .0154 &  1.54  &   1.40   &  1,010     &   224  &       64,526   &    111.24  &       27.58  & 68      \\
17  &   .0135 &  1.55  &   1.45   &  1,770     &   612  &       129,302  &    558.46  &       41.02  & 112     \\
18  &   .0037 &  1.46  &   1.33   &  980      &   180  &       261,164  &    52.93   &       32.45  & 107    \\
19  &   .0078 &  1.54  &   1.50   &  4,125     &   2,242 &       520,163  &    824.24  &       52.35  & 168     \\
20  &   .0011 &  1.42  &   1.32   &  1,227     &   280  &       1,047,349 &    88.00   &       77.69  & 196     \\
\hline 
21  &   .0008 &  1.42  &   1.39   &  1,731     &   1,092 &       2,095,421 &    29.02   &       180.81 & 480     \\
22  &   .0006 &  1.43  &   1.27   &  2,829     &   220  &       4,191,475 &    85.05   &       134.13 & 399     \\
23  &   .0008 &  1.46  &   1.44   &  6,833     &   4,462 &       8,381,775 &    4,148.57 &       209.22 & 699     \\
\hline
\end{tabular}
\vskip .1in
%this map was from the output of program inG to file G.txt 
\caption{The c.a. is $G=x_{-1}+x_0x_1+x_2$, which is 
bipermutative but not linear.}\label{ablebipermut}
\end{center}
\end{table}

\begin{table}[ht]
\begin{center}
\begin{tabular}{|r|r|r|r| r|r|r|r| r|r|}
\hline 
    &  Fraction  &  &          &  & &        &     Average   & Average     &  Maximum  \\
$k$ &  Periodic  & $\nu_k$       &  $L^{1/k}$     & P         &L & Not-P   &    Period     &  Preperiod     & Preperiod \\ 
\hline 
1   &  0.500000  & 1.00   &  1.00    & 1         &	1         &  1         &  1.00         &  0.50         &  1.00     \\
2   &  0.750000  & 1.73   &  1.00    & 3         &	1         &  1         &  1.00         &  0.25         &  1.00     \\
3   &  0.500000  & 1.58   &  1.00    & 4         &	1         &  4         &  1.00         &  0.88         &  2.00     \\
4   &  0.687500  & 1.82   &  1.41    & 11        &	4         &  5         &  2.50         &  0.31         &  1.00     \\
5   &  0.812500  & 1.91   &  1.71    & 26        &	15        &  6    &  9.75         &  0.19         &  1.00     	   \\
\hline		       		 										   
6   &  0.281250  & 1.61   &  1.00    & 18        &	1         &  46        &  1.00         &  3.95         &  8.00     \\
7   &  0.609375  & 1.86   &  1.74    & 78        &	49        &  50        &  37.75        &  0.94         &  5.00     \\
8   &  0.667969  & 1.90   &  1.81    & 171       &	120       &  85        &  75.94        &  0.86         &  6.00     \\
9   &  0.482422  & 1.84   &  1.55    & 247       &	54        &  265       &  44.93        &  2.83         &  12.00    \\
10  &  0.535156  & 1.87   &  1.82    & 548       &	410       &  476       &  345.04       &  2.85         &  17.00    \\
\hline		       		 										   
11  &  0.183105  & 1.71   &  1.60    & 375       &	176       &  1673      &  158.91       &  28.00        &  73.00    \\
12  &  0.176270  & 1.73   &  1.40    & 722       &	60        &  3374      &  6.38         &  37.95        &  85.00    \\
13  &  0.200073  & 1.76   &  1.59    & 1639      &	416       &  6553      &  220.46       &  19.73        &  76.00    \\
14  &  0.212524  & 1.79   &  1.62    & 3482      &	882       &  12902     &  483.97       &  42.97        &  153.00   \\
15  &  0.231598  & 1.81   &  1.59    & 7589      &	1095      &  25179     &  523.90       &  42.69        &  191.00   \\
\hline		       		 										   
16  &  0.117599  & 1.74   &  1.63    & 7707      &	2688      &  57829     &  1422.26      &  159.56       &  457.00   \\
17  &  0.078995  & 1.72   &  1.60    & 10354     &	3230      &  120718    &  2481.50      &  371.77       &  938.00   \\
18  &  0.078449  & 1.73   &  1.37    & 20565     &	324       &  241579    &  302.81       &  350.15       &  1155.00  \\
19  &  0.061646  & 1.72   &  1.64    & 32320     &	13471     &  491968    &  12128.71     &  404.87       &  1233.00  \\
20  &  0.065800  & 1.74   &  1.64    & 68996     &	21240     &  979580    &  15870.41     &  285.87       &  1063.00  \\
\hline		       		 										   
21  &  0.032823  & 1.69   &  1.56    & 68835     &	11865     &  2028317   &  816.87       &  1050.92      &  3506.00  \\
22  &  0.021364  & 1.67   &  1.60    & 89609     &	32428     &  4104695   &  20280.02     &  1335.34      &  5030.00  \\
23  &  0.011244  & 1.64   &  1.48    & 94324     &	9108      &  8294284   &  7929.18      &  4869.70      &  10024.00 \\
\hline 
\end{tabular}
\vskip .1in
%this map was from the output of program inB
\caption{The c.a. is $B=x_0+x_1x_2$ on the 2-shift. 
$B$  is degree 1, left permutative and
not right closing. }\label{ableleftpermmap}
\end{center}

\end{table}

\begin{table}[ht]
\begin{center}
\begin{tabular}{|r|r|r|r| r|r|r|r| r|r|}
\hline 
    &  Fraction  &  &          &  & &        &   Average   & Average     &  Maximum  \\
$k$ &  Periodic  & $\nu_k$          &  $L^{1/k}$  & P         &L & Not-P   &    Period     &  Preperiod     & Preperiod \\ 
\hline 
1   &  .5000 &  1.00  &   1.00   &  1        &   1      &     1       &    1.00       &    0.50        &   1     \\
2   &  .7500 &  1.73  &   1.00   &  3        &   1      &     1       &    1.00       &    0.25        &   1      \\
3   &  .5000 &  1.58  &   1.00   &  4        &   1      &     4       &    1.00       &    0.88        &   2      \\
4   &  .6875 &  1.82  &   1.00   &  11       &   1      &     5       &    1.00       &    0.31        &   1      \\
5   &  .8125 &  1.91  &   1.37   &  26       &   5      &     6       &    2.56       &    0.19        &   1     \\
\hline	    	    	      
6   &  .6562 &  1.86  &   1.61   &  42       &   18     &     22      &    7.38       &    0.86        &   4     \\
7   &  .6093 &  1.86  &   1.60   &  78       &   28     &     50      &    17.35      &    1.16        &   5      \\
8   &  .5117 &  1.83  &   1.62   &  131      &   48     &     125     &    27.50      &    2.08        &   10     \\
9   &  .4296 &  1.82  &   1.58   &  220      &   63     &     292     &    43.61      &    3.39        &   11     \\
10  &  .4082 &  1.82  &   1.31   &  418      &   15     &     606     &    5.45       &    6.91        &   21     \\
\hline	    	    	      
11  &  .4355 &  1.85  &   1.57   &  892      &   143    &     1,156    &    99.90      &    12.91       &   53     \\
12  &  .3608 &  1.83  &   1.36   &  1,478     &   42     &     2,618    &    11.61      &    15.59       &   53    \\
13  &  .3270 &  1.83  &   1.66   &  2,679     &   754    &     5,513    &    577.86     &    33.42       &   123    \\
14  &  .2167 &  1.79  &   1.26   &  3,552     &   28     &     12,832   &    23.06      &    79.16       &   191    \\
15  &  .2503 &  1.82  &   1.48   &  8,204     &   385    &     24,564   &    303.28     &    69.75       &   232    \\
\hline	    	    	      
16  &  .3152 &  1.86  &   1.63   &  20,659    &   2,528   &     44,877   &    1,197.54    &    48.40       &   281    \\
17  &  .1784 &  1.80  &   1.55   &  23,393    &   1,853   &     107,679  &    1,538.93    &    168.75      &   464    \\
18  &  .1821 &  1.81  &   1.59   &  47,760    &   4,464   &     214,384  &    3,208.77    &    172.00      &   697   \\
19  &  .1357 &  1.80  &   1.49   &  71,175    &   1,957   &     453,113  &    1,685.66    &    352.58      &   1082   \\
20  &  .1620 &  1.82  &   1.56   &  169,886   &   7,976   &     878,690  &    5,604.39    &    258.96      &   953    \\
\hline	    	    	      
21  &  .1032 &  1.79  &   1.52   &  216,612   &   7,056   &     1,880,540 &    6,344.22    &    2,389.64     &   4,363   \\
22  &  .0902 &  1.79  &   1.35   &  378,612   &   740    &     3,815,692 &    633.16     &    2,315.42     &   6,465   \\
23  &  .0858 &  1.79  &   1.58   &  720,246   &   39,353  &     7,668,362 &    36,059.28   &    1,760.56     &   5,984   \\
\hline
\end{tabular}
\vskip .1in
%this map was from the output of program inK into file B#U.txt
\caption{The c.a. is $K$, the composition which is the automorphism $U$ of 
Table \ref{ableF} followed by the left permutative, not right-closing
map $x_0 + x_1x_2$. The c.a. $K$ is left closing, not
right-closing,
and not permutative on either side.}\label{ableK}
\end{center}
\end{table}

\begin{table}[ht]
\begin{center}
\begin{tabular}{|r|r|r|r| r|r|r|r| r|r|}
\hline 
    &  Fraction  &  &          &  & &        &     Average   & Average     &  Maximum  \\
$k$ &  Periodic  &  $\nu_k$   &  $L^{1/k}$    & P         &L & Not-P   &    Period     &  Preperiod     & Preperiod \\ 
\hline 
1   &  0.500  & 1.00   &  1.00    & 1         &	1         &  1         &  1.00         &  0.50       & 1      \\
2   &  0.750  & 1.73   &  1.00    & 3         &	1         &  1         &  1.00         &  0.25       & 1      \\
3   &  0.500  & 1.58   &  1.44    & 4         &	3         &  4         &  1.75         &  0.50       & 1      \\
4   &  0.687  & 1.82   &  1.18    & 11        &	2         &  5         &  1.75         &  0.31       & 1      \\
5   &  0.812  & 1.91   &  1.71    & 26        &	15        &  6         &  11.00        &  0.19       & 1     \\
\hline 
6   &  0.468  & 1.76   &  1.20    & 30        &	3         &  34        &  1.66         &  1.28       & 3     \\
7   &  0.500  & 1.81   &  1.66    & 64        &	35        &  64        &  28.34        &  0.99       & 4      \\
8   &  0.667  & 1.90   &  1.63    & 171       &	52        &  85        &  30.39        &  0.52       & 3      \\
9   &  0.306  & 1.75   &  1.27    & 157       &	9         &  355       &  8.89         &  2.35       & 8      \\
10  &  0.261  & 1.74   &  1.49    & 268       &	55        &  756       &  20.18        &  6.75       & 18     \\
\hline 
11  &  0.387  & 1.83   &  1.57    & 793       &	143       &  1255      &  53.53        &  3.00       & 13     \\
12  &  0.088  & 1.63   &  1.16    & 362       &	6         &  3734      &  1.39         &  20.61      & 48    \\
13  &  0.150  & 1.72   &  1.63    & 1236      &	611       &  6956      &  259.75       &  20.15      & 78     \\
14  &  0.126  & 1.72   &  1.51    & 2068      &	329       &  14316     &  119.61       &  33.22      & 132    \\
15  &  0.091  & 1.70   &  1.50    & 3014      &	465       &  29754     &  414.94       &  44.45      & 138    \\
\hline 
16  &  0.092  & 1.72   &  1.50    & 6043      &	728       &  59493     &  650.33       &  101.66     & 282    \\
17  &  0.107  & 1.75   &  1.68    & 14145     &	6783      &  116927    &  3918.82      &  48.16      & 196    \\
18  &  0.060  & 1.71   &  1.58    & 15753     &	4095      &  246391    &  3406.78      &  110.78     & 396   \\
19  &  0.072  & 1.74   &  1.60    & 38191     &	7619      &  486097    &  6336.19      &  142.98     & 406    \\
20  &  0.038  & 1.69   &  1.54    & 40396     &	5780      &  1008180   &  1691.96      &  279.69     & 780    \\
\hline 
21  &  0.018  & 1.65   &  1.48    & 37867     &	4011      &  2059285   &  3961.81      &  705.45     & 1777   \\
22  &  0.017  & 1.66   &  1.51    & 75309     &	9658      &  4118995   &  4527.64      &  605.57     & 1770   \\
23  &  0.017  & 1.67   &  1.57    & 144096    &	34477     &  8244512   &  26857.88     &  1191.56    & 2687   \\
\hline
\end{tabular}
\vskip .1in
%this map was from the output of program inC
\caption{The c.a. is $C$, the composition $x_0x_1+x_2$ followed by $x_0+x_1x_2$, on the
2-shift. The c.a. $C$  is neither left nor 
 right closing. It is not known whether the periodic points of $C$ are
dense. }\label{ablenotclosing}
\end{center}
\end{table}

\clearpage

\begin{table}[ht]
\begin{center}
\begin{tabular}{|r|r|r|r| r|r|r|r| r|r|}
\hline 
    &  Fraction  &          &          &  & &        &     Average   & Average     &  Maximum  \\
$k$ &  Periodic  &  $\nu_k$ &  $L^{1/k}$     & P         &L & Not-P   &    Period     &  Preperiod     & Preperiod \\ 
\hline 
1   & 1.00 & 3.00    &2.00     &3          & 2            & 0        &  1.67        &  0.00    &   0.00      \\
2   & 1.00 & 3.00    &2.44     &9          & 6            & 0        &  4.56        &  0.00    &   0.00     \\
3   & 1.00 & 3.00    &2.46     &27         & 15           & 0        &  9.30        &  0.00    &   0.00     \\
4   & 0.11 & 1.73    &1.56     &9          & 6            & 72       &  4.56        &  0.89    &   1.00     \\
5   & 1.00 & 3.00    &2.77     &243        & 165          & 0        &  121.13      &  0.00    &   0.00     \\
6   & 1.00 & 3.00    &2.80     &729        & 486          & 0        &  334.54      &  0.00    &   0.00     \\
7   & 1.00 & 3.00    &2.57     &2187       & 742          & 0         & 401.62       & 0.00     &  0.00      \\
8   & 0.01 & 1.73    &1.58     &81         & 40           & 6480      & 33.50        & 4.24     &  9.00      \\
9   & 1.00 & 3.00    &2.24     &19683       &1469          &0          &1185.85       &0.00      & 0.00       \\
10  & 1.00 & 3.00    &2.76     &59049       &25865         &0          &22737.63      &0.00      & 0.00       \\
11  & 1.00 & 3.00    &2.92     &177147      &131857        &0          &109208.21     &0.00      & 0.00        \\
12  & 0.00 & 1.76    &1.67     &909         &486           &530532     &239.27        &45.51     & 133.00     \\
13  & 1.00 & 3.00    &2.79     &1594323     &631605        &0         & 291222.95     &0.00      & 0.00         \\
\hline
\end{tabular}
\vskip .1in
%this map was from the output of program inJ3 into the file A#W.txt
\caption{This map on the 3-shift is an automorphism $W$ followed 
by the degree 9 linear map $x_0+x_2$, where 
$W  = x_0+ 2x_0x_1x_1+ 2x_0x_1+x_1*x_1+x_1$. 
Let $\pi$ denote the permutation on $\{0,1,2\}$ which transposes 
$0$ and $2$. Then $(Wx)_0=x_0$ if $x_1\neq 1$ and $(Wx)_0=
\pi (x_0)$ if $x_1=1$.}   
\label{3shiftbiclosdeg9}
\end{center}
\end{table}

\begin{table}[ht]
\begin{center}
\begin{tabular}{|r|r|r|r| r|r|r|r| r|r|}
\hline 
    &  Fraction  &        &          &  & &        &     Average   & Average     &  Maximum  \\
$k$ &  Periodic  & $\nu_k$&  $L^{1/k}$     & P         &L & Not-P   &    Period     &  Preperiod     & Preperiod \\ 
\hline 
1   & .3333 &  1.00   &  1.00   &  1       &    1        &   2        &   1.00       &    1.00     &      2       \\
2   & .5555 &  2.23   &  1.00   &  5       &    1        &   4        &   1.00       &    0.56     &      2      \\
3   & .2592 &  1.91   &  1.00   &  7       &    1        &   20       &   1.00       &    1.78     &      3      \\
4   & .2098 &  2.03   &  1.00   &  17      &    1        &   64       &   1.00       &    3.17     &      8      \\
5   & .4362 &  2.54   &  2.09   &  106     &    40       &   137      &   27.65      &    1.08     &      4      \\
6   & .2208 &  2.33   &  1.51   &  161     &    12       &   568      &   6.79       &    2.71     &      9      \\
7   & .0932 &  2.13   &  1.66   &  204     &    35       &   1,983     &   5.89       &    13.89    &      38     \\
8   & .0391 &  2.00   &  1.00   &  257     &    1        &   6,304     &   1.00       &    27.02    &      67     \\
9   & .1667 &  2.45   &  1.60   &  3,283    &    72       &   16,400    &   48.89      &    13.41    &      52     \\
10  & .0299 &  2.11   &  1.62   &  1,770    &    130      &   57,279    &   62.67      &    55.38    &      163    \\
11  & .0224 &  2.12   &  1.89   &  3,972    &    1122     &   173,175   &   593.34     &    99.23    &      297     \\
12  & .0164 &  2.13   &  1.40   &  8,729    &    60       &   522,712   &   12.56      &    88.45    &      222    \\
13  & .0076 &  2.06   &  1.81   &  12,117   &    2366     &   1,582,206  &   2,228.50    &    676.85   &      1,504    \\
\hline
\end{tabular}
\vskip .1in
%this map was from the output of program inB3 into the file B3.txt
\caption{The map $x_0 + x_1x_2$ on the 3-shift: still degree 1, left
permutative, not right closing.}\label{3shift}
\end{center}
\end{table}

\begin{table}[ht]
\begin{center}
\begin{tabular}{|c|c|c|c|c|c|c|c|c|c|c|c|c|c|}
\hline
k     & 1&2&3&4&5&6&7&8&9&10&11&12&13\\ \hline 
1     &1.00     &1.00   &2.00  &2.00   &1.00   &2.00   &2.00   &2.00   &1.00   &1.00   &1.00   &1.00   &1.00\\
2     &0.00     &0.00   &0.00  &1.41   &1.41   &0.00   &0.00   &0.00   &0.00   &0.00   &0.00   &0.00   &0.00\\
3     &1.81     &1.44   &1.81  &1.81   &0.00   &1.44   &1.44   &1.81   &1.81   &1.44   &1.44   &0.00   &0.00\\
4     &1.41     &0.00   &0.00  &1.86   &0.00   &1.41   &1.68   &0.00   &0.00   &0.00   &1.68   &1.68   &1.41\\
5     &1.90     &1.82   &1.90  &1.82   &1.71   &1.97   &1.71   &1.97  &1.90   &1.37   &1.71   &1.71   &1.82\\ 
\hline 
6     &1.61     &0.00   &1.51  &1.69   &0.00   &1.69   &1.76   &1.86   &1.81   &1.61   &1.51   &1.61   &0.00\\
7     &1.77     &1.80   &1.80  &1.77   &1.74   &1.88   &1.90   &1.99   &1.45   &1.85   &1.90   &1.70   &1.92\\
8     &1.62     &1.70   &1.75  &1.70   &1.48   &1.72   &1.81   &1.90   &1.80   &0.00   &1.78   &1.72   &1.70\\
9     &1.89     &1.76   &1.74  &1.90   &1.90   &1.91   &1.91   &1.99   &1.79   &1.58   &1.76   &1.77   &1.80\\
10    &1.80     &1.58   &1.25  &1.83   &1.78   &1.81   &1.91   &1.95   &1.67   &0.00   &1.76   &1.68  &1.90\\
\hline
11    &1.66     &1.69   &1.75  &1.73   &1.78   &1.85   &1.92   &1.99   &1.67   &1.66   &1.80   &1.51   &1.48\\
12    &1.71     &1.75   &1.78  &1.84   &1.68   &1.85   &1.84   &1.97   &1.77   &1.70   &1.68   &1.73   &1.69\\
13    &1.73     &1.72   &1.79  &1.73   &1.72   &1.87   &1.93   &2.00   &1.72   &1.69   &1.80   &1.75   &1.84\\
14    &1.66     &1.61   &1.73  &1.73   &1.63   &1.81   &1.91   &1.98   &1.57   &1.54   &1.69   &1.68   &1.74\\
15    &1.66     &1.71   &1.60  &1.73   &1.74   &1.85   &1.92   &1.99   &1.78   &1.67   &1.70   &1.68   &1.76\\
\hline
16    &1.68     &1.64   &1.74  &1.71   &1.72   &1.79   &1.93   &1.98   &1.54   &1.65   &1.67   &1.49   &1.75\\
17    &1.69     &1.53   &1.73  &1.68   &1.72   &1.84   &1.91   &2.00   &1.73   &1.59   &1.73   &1.65   &1.69\\
18    &1.68     &1.46   &1.69  &1.68   &1.71   &1.83   &1.91   &1.99   &1.65   &1.59   &1.57   &1.65   &1.67\\
19    &1.67     &1.61   &1.68  &1.67   &1.69   &1.81   &1.93   &2.00   &1.71   &1.63   &1.66   &1.71   &1.73\\
\hline
\end{tabular}
\vskip .1in
\caption{$\nu^{o}_k(\cdot ,S_2)$ for the span five maps 1-13 of Table
\ref{span5maps}.}
\label{ablespan5a}
\end{center}
\end{table}

\begin{table}[ht]
\begin{center}
\begin{tabular}{|c|c|c| c|c|c| c|c|c|c| c|c|c|c|}
\hline
k     & 14& 15& 16& 17& 18& 19& 20& 21& 22& 23& 24& 25& 26\\
\hline
1  &1.00    &2.00  &2.00  &2.00  &1.00  &1.00  &1.00  &1.00  &1.00  &1.00  &2.00  &2.00  &2.00\\
2  &0.00    &1.41  &1.41  &0.00  &1.41  &1.41  &0.00  &0.00  &0.00  &0.00  &0.00  &0.00  &0.00\\
3  &1.81    &1.44  &0.00  &1.44  &1.44  &1.44  &0.00  &0.00  &0.00  &0.00  &1.81  &1.81  &1.81\\
4  &0.00    &1.68  &1.86  &0.00  &1.41  &1.68  &0.00  &1.41  &0.00  &1.68  &0.00  &0.00  &0.00\\
5  &1.82    &1.82  &1.82  &1.71  &1.71  &1.90  &1.71  &1.58  &1.71  &1.58  &1.97  &1.97  &1.97\\
\hline
6  &1.81    &1.86  &1.76  &1.61  &1.51  &1.51  &0.00  &1.69  &1.81  &1.81  &1.69  &1.69  &1.69\\
7  &1.80    &1.85  &1.92  &1.60  &1.70  &1.83  &1.83  &1.54  &1.70  &1.54  &1.99  &1.99  &1.99\\
8  &1.41    &1.83  &1.70  &1.72  &1.83  &1.86  &1.68  &1.68  &1.54  &1.68  &1.81  &1.81  &1.86\\
9  &1.84    &1.84  &1.76  &1.76  &1.71  &1.78  &1.78  &1.74  &1.84  &1.87  &1.99  &1.99  &1.99\\
10  &1.86    &1.85  &1.83  &1.65  &1.76  &1.82  &1.79  &1.52  &1.67  &1.72  &1.90  &1.90  &1.93\\
\hline
11  &1.57    &1.78  &1.81  &1.88  &1.81  &1.75  &1.77  &1.64  &1.79  &1.54  &1.99  &1.99  &1.99\\
12  &1.57    &1.83  &1.84  &1.66  &1.64  &1.80  &1.62  &1.73  &1.73  &1.67  &1.94  &1.94  &1.94\\
13  &1.71    &1.71  &1.68  &1.73  &1.72  &1.85  &1.61  &1.78  &1.79  &1.65  &2.00  &2.00  &2.00\\
14  &1.76    &1.73  &1.72  &1.67  &1.66  &1.80  &1.57  &1.73  &1.75  &1.60  &1.95  &1.95  &1.96\\
15  &1.71    &1.76  &1.79  &1.72  &1.77  &1.78  &1.72  &1.65  &1.66  &1.52  &1.99  &1.99  &1.99\\
\hline
16  &1.71    &1.71  &1.71  &1.73  &1.69  &1.78  &1.71  &1.51  &1.65  &1.55  &1.97  &1.97  &1.97\\
17  &1.74    &1.66  &1.63  &1.64  &1.71  &1.78  &1.64  &1.65  &1.67  &1.68  &2.00  &2.00  &2.00\\
18  &1.71    &1.68  &1.68  &1.67  &1.62  &1.73  &1.65  &1.61  &1.52  &1.62  &1.98  &1.98  &1.98\\
19  &1.68    &1.69  &1.70  &1.70  &1.63  &1.75  &1.70  &1.66  &1.65  &1.60  &2.00  &2.00  &2.00\\
\hline
\end{tabular}
\vskip .1in
\caption{$\nu^{o}_k(\cdot ,S_2)$ for the span five  maps 14-26 of
Table \ref{span5maps}.}\label{ablespan5b}
\end{center}
\end{table}

\begin{table}[ht]
\begin{center}
\begin{tabular}{|c|c|c|c|c |c|c|c|c|c |c|c|c|c|c |c|c|}
\hline
k     &1 &2 &3 &4 &5 &6 &7 &8 &9 &10 &11 &12 &13 &14 &15 &16\\
\hline
1     &2.00&2.00 &1.00  &1.00  &2.00  &2.00  &2.00  &1.00  &2.00  &1.00  &1.00  &2.00  &2.00  &1.00  &1.00  & 2.00   \\
2     &0.00&1.41 &0.00  &1.41  &1.41  &1.41  &1.41  &1.41  &0.00  &0.00  &0.00  &1.41  &0.00  &1.41  &1.41  & 1.41   \\
3     &1.81&1.81 &1.81  &1.44  &0.00  &1.81  &1.81  &1.44  &1.81  &1.44  &1.81  &0.00  &1.81  &1.44  &1.44  & 1.81   \\
4     &1.68&1.86 &1.68  &1.41  &1.41  &1.86  &1.86  &1.41  &1.68  &1.41  &1.68  &1.41  &1.68  &0.00  &1.68  & 1.86   \\
5     &1.97&1.97 &1.82  &1.90  &1.58  &1.97  &1.97  &1.90  &1.97  &1.82  &1.82  &1.58  &1.97  &1.58  &1.90  & 1.97   \\
\hline
6     &1.90&1.94 &1.81  &1.61  &1.76  &1.94  &1.94  &1.61  &1.90  &1.34  &1.81  &1.76  &1.90  &1.86  &1.81  & 1.94   \\
7     &1.99&1.99 &1.83  &1.70  &1.80  &1.99  &1.99  &1.70  &1.99  &1.74  &1.83  &1.80  &1.99  &1.83  &1.70  & 1.99   \\
8     &1.95&1.98 &1.81  &1.48  &1.75  &1.98  &1.98  &1.48  &1.95  &1.78  &1.81  &1.75  &1.95  &1.54  &1.88  & 1.98   \\
9     &1.99&1.99 &1.86  &1.68  &1.82  &1.99  &1.99  &1.68  &1.99  &1.82  &1.86  &1.82  &1.99  &1.73  &1.86  & 1.99   \\
10    &1.98&1.99 &1.76  &1.70  &1.82  &1.99  &1.99  &1.70  &1.98  &1.75  &1.76  &1.82  &1.98  &1.56  &1.65  & 1.99   \\
\hline
11    &1.99&1.99 &1.70  &1.65  &1.68  &1.99  &1.99  &1.65  &1.99  &1.89  &1.70  &1.68  &1.99  &1.60  &1.90  & 1.99   \\
12    &1.99&1.99 &1.51  &1.65  &1.61  &1.99  &1.99  &1.65  &1.99  &1.65  &1.51  &1.61  &1.99  &1.34  &1.75  & 1.99   \\
13    &2.00&2.00 &1.70  &1.57  &1.63  &2.00  &2.00  &1.57  &2.00  &1.73  &1.70  &1.63  &2.00  &1.54  &1.68  & 2.00   \\
14    &1.99&1.99 &1.74  &1.65  &1.70  &1.99  &1.99  &1.65  &1.99  &1.81  &1.74  &1.70  &1.99  &1.66  &1.74  & 1.99   \\
15    &1.99&1.99 &1.71  &1.68  &1.70  &1.99  &1.99  &1.68  &1.99  &1.73  &1.71  &1.70  &1.99  &1.47  &1.77  & 1.99   \\
\hline
16    &1.99&1.99 &1.74  &1.67  &1.70  &1.99  &1.99  &1.67  &1.99  &1.76  &1.74  &1.70  &1.99  &1.59  &1.67  & 1.99   \\
17    &2.00&2.00 &1.67  &1.53  &1.71  &2.00  &2.00  &1.53  &2.00  &1.75  &1.67  &1.71  &2.00  &1.59  &1.61  & 2.00   \\
18    &1.99&1.99 &1.71  &1.56  &1.65  &1.99  &1.99  &1.56  &1.99  &1.71  &1.71  &1.65  &1.99  &1.52  &1.63  & 1.99   \\
19    &2.00&2.00 &1.73  &1.54  &1.72  &2.00  &2.00  &1.54  &2.00  &1.77  &1.73  &1.72  &2.00  &1.57  &1.69  & 2.00   \\
\hline
\end{tabular}
\vskip .1in 
\caption{$\nu^{o}_k(\cdot ,S_2)$ for the span four maps 1-16
of Table \ref{span4maps}.}\label{ablespan4a}
\end{center}
\end{table}

\begin{table}[ht]
\begin{center}
\begin{tabular}{|c|c|c|c|c |c|c|c|c|c |c|c|c|c|c |c|c|}
\hline
k     &$F1$ &$F2$  &$F3$  &$F4$  &$F5$  &$F6$  &$F7$  &$F8$ 
      &$F9$ &$F10$ &$F11$ &$F12$ &$F13$ &$F14$ &$F15$ &$F16$\\

\hline
1     &2.00   &2.00  &1.00  &1.00  &2.00  &2.00 &2.00  &1.00  &2.00  &1.00  &1.00  &2.00  &2.00  &1.00 & 1.00 &2.00   \\   
2     &0.00   &1.41  &0.00  &1.41  &1.41  &1.41 &1.41  &1.41  &0.00  &0.00  &0.00  &1.41  &0.00  &1.41 & 1.41 &1.41   \\   
3     &1.81   &1.81  &1.81  &1.44  &1.44  &1.81 &1.81  &1.44  &1.81  &0.00  &1.81  &1.44  &1.81  &1.44 & 0.00 &1.81   \\   
4     &1.41   &1.86  &1.68  &0.00  &1.41  &1.86 &1.86  &0.00  &1.41  &0.00  &1.68  &1.41  &1.41  &1.41 & 1.68 &1.86   \\   
5     &1.97   &1.97  &1.82  &1.58  &1.71  &1.97 &1.97  &1.58  &1.97  &1.71  &1.82  &1.71  &1.97  &1.90 & 1.90 &1.97   \\   
\hline
6     &1.86   &1.94  &1.69  &1.86  &1.69  &1.94 &1.94  &1.86  &1.86  &1.69  &1.69  &1.69  &1.86  &1.61 & 0.00 &1.94   \\   
7     &1.99   &1.99  &1.66  &1.83  &1.70  &1.99 &1.99  &1.83  &1.99  &1.60  &1.66  &1.70  &1.99  &1.70 & 1.92 &1.99   \\   
8     &1.93   &1.98  &1.81  &1.54  &1.80  &1.98 &1.98  &1.54  &1.93  &1.41  &1.81  &1.80  &1.93  &1.48 & 1.75 &1.98   \\   
9     &1.99   &1.99  &1.71  &1.73  &1.62  &1.99 &1.99  &1.73  &1.99  &1.77  &1.71  &1.62  &1.99  &1.68 & 1.68 &1.99   \\   
10    &1.95   &1.99  &1.70  &1.56  &1.44  &1.99 &1.99  &1.56  &1.95  &1.79  &1.70  &1.44  &1.95  &1.70 & 0.00 &1.99   \\   
\hline
11    &1.99   &1.99  &1.74  &1.60  &1.46  &1.99 &1.99  &1.60  &1.99  &1.59  &1.74  &1.46  &1.99  &1.65 & 1.65 &1.99   \\   
12    &1.97   &1.99  &1.65  &1.34  &1.55  &1.99 &1.99  &1.34  &1.97  &1.63  &1.65  &1.55  &1.97  &1.65 & 1.69 &1.99   \\   
13    &2.00   &2.00  &1.75  &1.54  &1.65  &2.00 &2.00  &1.54  &2.00  &1.53  &1.75  &1.65  &2.00  &1.57 & 1.67 &2.00   \\   
14    &1.98   &1.99  &1.74  &1.66  &1.53  &1.99 &1.99  &1.66  &1.98  &1.72  &1.74  &1.53  &1.98  &1.65 & 1.51 &1.99   \\   
15    &1.99   &1.99  &1.74  &1.47  &1.64  &1.99 &1.99  &1.47  &1.99  &1.68  &1.74  &1.64  &1.99  &1.68 & 1.74 &1.99   \\   
\hline
16    &1.99   &1.99  &1.66  &1.59  &1.55  &1.99 &1.99  &1.59  &1.99  &1.66  &1.66  &1.55  &1.99  &1.67 & 1.68 &1.99   \\   
17    &2.00   &2.00  &1.67  &1.59  &1.57  &2.00 &2.00  &1.59  &2.00  &1.74  &1.67  &1.57  &2.00  &1.53 & 1.69 &2.00   \\   
18    &1.99   &1.99  &1.63  &1.52  &1.61  &1.99 &1.99  &1.52  &1.99  &1.70  &1.63  &1.61  &1.99  &1.56 & 1.61 &1.99   \\   
19    &2.00   &2.00  &1.69  &1.57  &1.51  &2.00 &2.00  &1.57  &2.00  &1.54  &1.69  &1.51  &2.00  &1.54 & 1.63 &2.00    \\     
\hline
\end{tabular}
\caption{$\nu^{o}_k(\cdot ,S_2)$ for the span 4 maps 1-16 of Table
\ref{span4maps}, postcomposed with 
the flip map $F=x_0+1$.}\label{ablespan4c} 
\end{center}
\end{table}

\begin{table}[ht]
\begin{center}
\begin{tabular}{|c|c|c|c|c |c|c|c|c|c |c|c|c|c|c |c|c|}
\hline
k     &17 &18 &19 &20 &21 &22 &23 &24 &25 &26 &27 &28 &29 &30 &31 &32\\
\hline 
1    &1.00  &1.00  &2.00  &2.00  &1.00   &2.00  &1.00  &1.00  &2.00  &2.00  &1.00  &1.00  &2.00  &2.00  &1.00  & 1.00  \\
2    &0.00  &1.41  &0.00  &1.41  &1.41   &0.00  &0.00  &1.41  &1.41  &0.00  &0.00  &1.41  &0.00  &0.00  &0.00  & 0.00  \\
3    &0.00  &0.00  &0.00  &0.00  &1.44   &0.00  &1.81  &0.00  &0.00  &1.81  &1.44  &1.44  &0.00  &0.00  &0.00  & 1.81  \\
4    &0.00  &1.68  &1.68  &1.41  &1.68   &1.68  &1.68  &1.68  &1.41  &1.68  &1.41  &0.00  &1.68  &1.68  &0.00  & 1.68  \\
5    &1.71  &1.90  &1.71  &1.58  &1.90   &1.71  &1.82  &1.90  &1.58  &1.97  &1.82  &1.58  &1.71  &1.71  &1.71  & 1.82  \\
\hline 												  
6    &1.69  &0.00  &1.34  &1.76  &1.81   &1.34  &1.69  &0.00  &1.76  &1.90  &1.34  &1.86  &1.34  &1.34  &1.69  & 1.69  \\
7    &1.60  &1.92  &1.54  &1.80  &1.70   &1.54  &1.66  &1.92  &1.80  &1.99  &1.74  &1.83  &1.54  &1.54  &1.60  & 1.66  \\
8    &1.41  &1.75  &1.41  &1.75  &1.88   &1.41  &1.81  &1.75  &1.75  &1.95  &1.78  &1.54  &1.41  &1.41  &1.41  & 1.81  \\
9    &1.77  &1.68  &1.74  &1.82  &1.86   &1.74  &1.71  &1.68  &1.82  &1.99  &1.82  &1.73  &1.74  &1.74  &1.77  & 1.71  \\
10   &1.79  &0.00  &1.72  &1.82  &1.65   &1.72  &1.70  &0.00  &1.82  &1.98  &1.75  &1.56  &1.72  &1.72  &1.79  & 1.70  \\
\hline 												  
11   &1.59  &1.65  &1.41  &1.68  &1.90   &1.41  &1.74  &1.65  &1.68  &1.99  &1.89  &1.60  &1.41  &1.41  &1.59  & 1.74  \\
12   &1.63  &1.69  &1.59  &1.61  &1.75   &1.59  &1.65  &1.69  &1.61  &1.99  &1.65  &1.34  &1.59  &1.59  &1.63  & 1.65  \\
13   &1.53  &1.67  &1.66  &1.63  &1.68   &1.66  &1.75  &1.67  &1.63  &2.00  &1.73  &1.54  &1.66  &1.66  &1.53  & 1.75  \\
14   &1.72  &1.51  &1.44  &1.70  &1.74   &1.44  &1.74  &1.51  &1.70  &1.99  &1.81  &1.66  &1.44  &1.44  &1.72  & 1.74  \\
15   &1.68  &1.74  &1.58  &1.70  &1.77   &1.58  &1.74  &1.74  &1.70  &1.99  &1.73  &1.47  &1.58  &1.58  &1.68  & 1.74  \\
\hline 												  
16   &1.66  &1.68  &1.64  &1.70  &1.67   &1.64  &1.66  &1.68  &1.70  &1.99  &1.76  &1.59  &1.64  &1.64  &1.66  & 1.66  \\
17   &1.74  &1.69  &1.59  &1.71  &1.61   &1.59  &1.67  &1.69  &1.71  &2.00  &1.75  &1.59  &1.59  &1.59  &1.74  & 1.67  \\
18   &1.70  &1.61  &1.46  &1.65  &1.63   &1.46  &1.63  &1.61  &1.65  &1.99  &1.71  &1.52  &1.46  &1.46  &1.70  & 1.63  \\
19   &1.54  &1.63  &1.60  &1.72  &1.69   &1.60  &1.69  &1.63  &1.72  &2.00  &1.77  &1.57  &1.60  &1.60  &1.54  & 1.69  \\

\hline
\end{tabular}
\vskip .1in 
\caption{$\nu^{o}_k(\cdot ,S_2)$ for the span 4 onto maps 17-32
of Table \ref{span4maps}.}\label{ablespan4b}
\end{center}
\end{table}

\begin{table}[ht]
\begin{center}
\begin{tabular}{|c|c|c|c|c |c|c|c|c|c |c|c|c|c|c |c|c|}
\hline
k    &$F17$ &$F18$ &$F19$ &$F20$ &$F21$ &$F22$ &$F23$ &$F24$ 
     &$F25$ &$F26$ &$F27$ &$F28$ &$F29$ &$F30$ &$F31$ &$F32$\\ 
\hline
1     &1.00  &1.00  &2.00  &2.00  &1.00  &2.00  &1.00  &1.00  &2.00  &2.00  &1.00  &1.00  &2.00  &2.00  & 1.00 &1.00  \\
2     &0.00  &1.41  &0.00  &1.41  &1.41  &0.00  &0.00  &1.41  &1.41  &0.00  &0.00  &1.41  &0.00  &0.00  & 0.00 &0.00  \\
3     &1.44  &1.44  &1.44  &1.44  &0.00  &1.44  &1.81  &1.44  &1.44  &1.81  &0.00  &1.44  &1.44  &1.44  & 1.44 &1.81  \\
4     &1.41  &1.68  &1.68  &1.41  &1.68  &1.68  &1.68  &1.68  &1.41  &1.41  &0.00  &1.41  &1.68  &1.68  & 1.41 &1.68  \\
5     &1.82  &1.90  &1.58  &1.71  &1.90  &1.58  &1.82  &1.90  &1.71  &1.97  &1.71  &1.90  &1.58  &1.58  & 1.82 &1.82  \\
\hline 												  
6     &1.34  &1.81  &1.61  &1.69  &0.00  &1.61  &1.81  &1.81  &1.69  &1.86  &1.69  &1.61  &1.61  &1.61  & 1.34 &1.81  \\
7     &1.74  &1.70  &1.70  &1.70  &1.92  &1.70  &1.83  &1.70  &1.70  &1.99  &1.60  &1.70  &1.70  &1.70  & 1.74 &1.83  \\
8     &1.78  &1.88  &1.65  &1.80  &1.75  &1.65  &1.81  &1.88  &1.80  &1.93  &1.41  &1.48  &1.65  &1.65  & 1.78 &1.81  \\
9     &1.82  &1.86  &1.60  &1.62  &1.68  &1.60  &1.86  &1.86  &1.62  &1.99  &1.77  &1.68  &1.60  &1.60  & 1.82 &1.86  \\
10    &1.75  &1.65  &1.61  &1.44  &0.00  &1.61  &1.76  &1.65  &1.44  &1.95  &1.79  &1.70  &1.61  &1.61  & 1.75 &1.76  \\
\hline 												  
11    &1.89  &1.90  &1.69  &1.46  &1.65  &1.69  &1.70  &1.90  &1.46  &1.99  &1.59  &1.65  &1.69  &1.69  & 1.89 &1.70  \\
12    &1.65  &1.75  &1.49  &1.55  &1.69  &1.49  &1.51  &1.75  &1.55  &1.97  &1.63  &1.65  &1.49  &1.49  & 1.65 &1.51  \\
13    &1.73  &1.68  &1.66  &1.65  &1.67  &1.66  &1.70  &1.68  &1.65  &2.00  &1.53  &1.57  &1.66  &1.66  & 1.73 &1.70  \\
14    &1.81  &1.74  &1.59  &1.53  &1.51  &1.59  &1.74  &1.74  &1.53  &1.98  &1.72  &1.65  &1.59  &1.59  & 1.81 &1.74  \\
15    &1.73  &1.77  &1.53  &1.64  &1.74  &1.53  &1.71  &1.77  &1.64  &1.99  &1.68  &1.68  &1.53  &1.53  & 1.73 &1.71  \\
\hline 												  
16    &1.76  &1.67  &1.69  &1.55  &1.68  &1.69  &1.74  &1.67  &1.55  &1.99  &1.66  &1.67  &1.69  &1.69  & 1.76 &1.74  \\
17    &1.75  &1.61  &1.63  &1.57  &1.69  &1.63  &1.67  &1.61  &1.57  &2.00  &1.74  &1.53  &1.63  &1.63  & 1.75 &1.67  \\
18    &1.71  &1.63  &1.55  &1.61  &1.61  &1.55  &1.71  &1.63  &1.61  &1.99  &1.70  &1.56  &1.55  &1.55  & 1.71 &1.71  \\
19    &1.77  &1.69  &1.65  &1.51  &1.63  &1.65  &1.73  &1.69  &1.51  &2.00  &1.54  &1.54  &1.65  &1.65  & 1.77 &1.73  \\
% the next row for k=20 is correct but causes space problems in the document 
% 20    &1.69  &1.65  &1.64  &1.54  &1.65  &1.64  &1.69  &1.65  &1.54  &1.99  &1.64  &1.58  &1.64  &1.56  & 1.69 &1.69  \\
\hline 
\end{tabular}
\vskip .1in
\caption{$\nu^{o}_k(\cdot ,S_2)$ for the span 4 maps 17-32 of Table 
\ref{span4maps}, postcomposed with 
the flip map $F=x_0+1$.}\label{ablespan4d}
\end{center}
\end{table}

\clearpage

\begin{table}[ht]
\begin{center}
\begin{tabular}{|c|c|c|c|c |c|c|c|c|c |c|c|c|c|c |c|c|}
\hline
k     &1 &2 &3 &4 &5 &6 &7 &8 &9 &10 &11 &12 &13 &14 &15 &16\\
\hline
1     &1.00 &1.00  &2.00  &2.00  &1.00  &1.00  &1.00  &2.00  &1.00  &2.00  &2.00  &1.00  &1.00  &2.00  &2.00  &1.00   \\
2     &1.41 &0.00  &1.41  &0.00  &0.00  &0.00  &0.00  &0.00  &1.41  &1.41  &1.41  &0.00  &1.41  &0.00  &0.00  &0.00   \\
3     &1.44 &0.00  &1.44  &0.00  &0.00  &1.44  &1.44  &1.81  &0.00  &0.00  &1.44  &1.81  &0.00  &1.81  &1.44  &1.44   \\
4     &0.00 &1.41  &1.41  &1.68  &1.68  &0.00  &1.86  &1.41  &0.00  &0.00  &1.41  &1.41  &1.41  &1.41  &1.68  &1.86   \\
5     &1.37 &1.37  &1.58  &1.71  &1.82  &1.90  &1.37  &1.71  &1.90  &1.71  &1.58  &1.82  &1.37  &1.71  &1.37  &1.90   \\
\hline
6     &1.51 &1.51  &1.51  &1.34  &1.51  &1.34  &1.34  &0.00  &0.00  &1.51  &0.00  &1.61  &1.34  &0.00  &1.34  &0.00   \\
7     &1.74 &1.88  &1.54  &1.60  &1.66  &1.60  &1.77  &1.54  &1.70  &1.66  &0.00  &1.74  &1.80  &1.54  &1.60  &1.32   \\
8     &1.48 &1.58  &1.65  &1.62  &1.41  &1.65  &0.00  &1.62  &1.62  &1.58  &1.68  &1.54  &0.00  &1.62  &1.68  &0.00   \\
9     &1.62 &1.78  &1.74  &1.58  &1.74  &1.58  &1.68  &0.00  &1.48  &1.48  &1.44  &1.73  &1.44  &0.00  &1.64  &1.60   \\
10    &1.34 &1.62  &1.40  &1.54  &1.66  &1.58  &1.52  &1.68  &1.61  &1.62  &1.34  &1.47  &1.34  &1.68  &1.34  &1.50   \\
\hline
11    &1.58 &1.69  &1.43  &1.64  &1.51  &1.76  &1.67  &1.60  &1.55  &1.53  &1.60  &1.66  &1.58  &1.60  &1.57  &1.67   \\
12    &1.54 &1.44  &1.60  &1.51  &1.63  &1.74  &1.46  &1.56  &0.00  &1.54  &1.54  &1.44  &1.52  &1.56  &1.38  &1.38   \\
13    &1.61 &1.66  &1.62  &1.65  &1.56  &1.77  &1.59  &1.63  &1.32  &1.66  &1.61  &1.47  &1.57  &1.63  &1.68  &1.60   \\
14    &1.64 &1.62  &1.52  &1.57  &1.53  &1.80  &1.55  &1.40  &1.55  &1.41  &1.60  &1.55  &1.48  &1.40  &1.55  &1.65   \\
15    &1.60 &1.74  &1.60  &1.38  &1.62  &1.73  &1.58  &1.62  &1.52  &1.52  &1.50  &1.49  &1.66  &1.62  &1.66  &1.61   \\
\hline
16    &1.48 &1.73  &1.47  &1.53  &1.50  &1.60  &1.51  &1.59  &1.49  &1.44  &1.29  &1.58  &1.57  &1.59  &1.49  &1.57   \\
17    &1.46 &1.69  &1.58  &1.59  &1.62  &1.64  &1.61  &1.47  &1.60  &1.58  &1.56  &1.55  &1.46  &1.47  &1.63  &1.36   \\
18    &1.56 &1.72  &1.55  &1.49  &1.55  &1.52  &1.46  &1.50  &1.51  &1.52  &1.54  &1.45  &1.56  &1.50  &1.52  &1.45   \\
19    &1.56 &1.66  &1.56  &1.55  &1.52  &1.64  &1.45  &1.60  &1.47  &1.51  &1.55  &1.56  &1.52  &1.60  &1.59  &1.48   \\
\hline
\end{tabular}
\vskip .1in
%This data was generated from input file inM into the directory ~/perca/resolv 
\caption{$\nu^o_k$ for 16 left permutative span 5 maps.
\newline 
Let $p_n(x[0,4])$ denote the polynomial rule  for 
the map $n$ 
in Table \ref{ableresolv} and let $q_n(x[0,3])$ denote the 
polynomial rule for the span 4 map $n$ in  
Table \ref{span4maps}. Then $p_n$ is defined by 
$p_n(x[0,4])=x_0+q_n(x[1,4])$.
\newline 
 The purpose of Table \ref{ableresolv}
is to allow a rough comparison of a sample of  maps which are 
linear in an end variable to maps which are not
(Table \ref{ableresolv} vs.
Table \ref{span4maps}). 
 We see no particular difference. 
}\label{ableresolv}
\end{center}
\end{table}

Table \ref{fpBcycles} gives complete cycle data for 
the c.a. $B$ through shift period $k=22$. 
of Table \ref{fppBsamplecycles}. 
For each $k$, 
all $B$-periods $p$  of points from 
$P_k(S_2)$ are listed. The multiplicities 
given are the number $\mu_{\textnormal{orb}}$ of $B$-cycles in 
$P_k(S_2)$ with the given size $p$; the number
$\mu_{\textnormal{per}}$ of points in all these cycles; 
and the number  $\mu_{\textnormal{ev}}$      
of points in $P_k(S_2)$ with eventual period $p$.

\begin{table}[ht]
\begin{center}
\begin{tabular}{|| r|r|r|r|r || r|r|r|r|r ||}
\hline  
k  & p & $\mu_{\textnormal{orb}}$&$\mu_{\textnormal{per}}$  &$\mu_{\textnormal{ev}}$       &
k  & p & $\mu_{\textnormal{orb}}$&$\mu_{\textnormal{per}}$  &  $\mu_{\textnormal{ev}}$      \\ 
\hline 								     								      
1  &   1       &      1             &  1               &   2&16 &   1       &      2207          &  2207            &   2208\\	
\cline{1-5}  
2  &   1       &      3             &  3               &   4         &	 &   4       &      1             &  4               &   6192      \\	  
\cline{1-5}  
3  &   1       &      4             &  4               &   8         &	 &   120     &      1             &  120             &   23520     \\	  
\cline{1-5}  
4  &   1       &      7             &  7               &   8         &	 &   2688    &      2             &  5376            &   33616     \\	  
\cline{6-10}
   &   4       &      1             &  4               &   8         &17 &   1       &      3571          &  3571            &   3572      \\	  
\cline{1-5}  
5  &   1       &      11            &  11              &   12        &	 &   1020    &      1             &  1020            &   1530      \\	  
   &   15      &      1             &  15              &   20        &	 &   2533    &      1             &  2533            &   119357    \\	  
\cline{1-5}  
6  &   1       &      18            &  18              &   64        &	 &   3230    &      1             &  3230            &   6613      \\	  
\cline{1-5}  \cline{6-10}
7  &   1       &      29            &  29              &   30        &18 &   1       &      5778          &  5778            &   5824      \\	  
   &   49      &      1             &  49              &   98        &	 &   9       &      1             &  9               &   9         \\	  
\cline{1-5}  
8  &   1       &      47            &  47              &   48        &	 &   38      &      36            &  1368            &   3834      \\	  
   &   4       &      1             &  4               &   48        &	 &   54      &      3             &  162             &   4815      \\	  
   &   120     &      1             &  120             &   160       &	 &   108     &      6             &  648             &   648       \\	  
\cline{1-5}  
9  &   1       &      76            &  76              &   80        &	 &   216     &      6             &  1296            &   7740      \\	  
   &   9       &      1             &  9               &   9         &	 &   296     &      36            &  10656           &   10656     \\	  
   &   54      &      3             &  162             &   423       &	 &   324     &      2             &  648             &   228618    \\	  
\cline{1-5}  \cline{6-10}
10 &   1       &      123           &  123             &   124       &19 &   1       &      9349          &  9349            &   9350      \\	  
   &   15      &      1             &  15              &   40        &	 &   76      &      1             &  76              &   76        \\	  
   &   410     &      1             &  410             &   860       &	 &   133     &      2             &  266             &   14421     \\	  
\cline{1-5}  
11 &   1       &      199           &  199             &   200       &	 &   171     &      1             &  171             &   171       \\	  
   &   176     &      1             &  176             &   1848      &	 &   646     &      1             &  646             &   2755      \\	  
\cline{1-5}  
12 &   1       &      322           &  322             &   3692      &	 &   4161    &      1             &  4161            &   25156     \\	  
   &   4       &      1             &  4               &   8         &	 &   4180    &      1             &  4180            &   12122     \\	  
   &   56      &      6             &  336             &   336       &	 &   13471   &      1             &  13471           &   460237    \\	  
\cline{6-10}
   &   60      &      1             &  60              &   60        &20 &   1       &      15127         &  15127           &   15128     \\	  
\cline{1-5}  
13 &   1       &      521           &  521             &   522       &	 &   4       &      1             &  4               &   8         \\	  
   &   10      &      13            &  130             &   650       &	 &   15      &      1             &  15              &   1000      \\	  
   &   26      &      1             &  26              &   117       &	 &   132     &      5             &  660             &   1560      \\	  
   &   143     &      1             &  143             &   3900      &	 &   140     &      4             &  560             &   560       \\	  
   &   403     &      1             &  403             &   845       &	 &   410     &      1             &  410             &   9420      \\	  
   &   416     &      1             &  416             &   2158      &	 &   5240    &      2             &  10480           &   306300    \\	  
\cline{1-5}  
14 &   1       &      843           &  843             &   844       &	 &   20500   &      1             &  20500           &   197240    \\	  
   &   49      &      1             &  49              &   602       &	 &   21240   &      1             &  21240           &   517360    \\	  
\cline{6-10}
   &   161     &      2             &  322             &   2212      &21 &   1       &      24476         &  24476           &   24480     \\	  
   &   448     &      2             &  896             &   7686      &	 &   14      &      3             &  42              &   42        \\	  
   &   490     &      1             &  490             &   882       &	 &   21      &      2             &  42              &   42        \\	  
   &   882     &      1             &  882             &   4158      &	 &   49      &      1             &  49              &   98        \\	  
\cline{1-5}  
15 &   1       &      1364          &  1364            &   1368      &	 &   57      &      21            &  1197            &   1197      \\	  
   &   15      &      1             &  15              &   20        &	 &   266     &      3             &  798             &   1949766   \\	  
   &   180     &      1             &  180             &   180       &	 &   2618    &      6             &  15708           &   15708     \\	  
   &   399     &      5             &  1995            &   8625      &	 &   4886    &      3             &  14658           &   14658     \\	  
   &   450     &      3             &  1350            &   15705     &	 &   11865   &      1             &  11865           &   91161     \\	  
\cline{6-10}
   &   530     &      3             &  1590            &   1590      &22 &   1       &      39603         &  39603           &   39604     \\	  
   &   1095    &      1             &  1095            &   5280      &	 &   132     &      2             &  264             &   264       \\	  
   &	       &		    &		       &             &	 &   176     &      1             &  176             &   910272    \\	  
   &	       &		    &		       &             &	 &   660     &      1             &  660             &   1100      \\	  
   &	       &		    &		       &	     &	 &   1067    &      2             &  2134            &   112948    \\	  
   &	       &		    &		       &	     &	 &   14344   &      1             &  14344           &   924814    \\	  
   &	       &		    &		       &	     &	 &   32428   &      1             &  32428           &   2205302   \\	  
\hline				      
\end{tabular}
\vskip .1in
\caption{
Table for $B=x_0+x_1x_2$ constructed from FPeriod: 
complete data  for shift-periods through $k=22$. See 
previous page for definitions. 
} 
\label{fpBcycles}
\end{center}
\end{table}

\clearpage

\section{FProbPeriod Tables} \label{fprobperiodtables}

For a given map $f$, a given positive integer $N$ 
and a given positive integer $m$ and given set of
positive integers $k$, 
the program FProbPeriod will for each $k$ randomly sample $m$ blocks
of length $k$ on alphabet $\{0, 1, \dots N-1\}$,
 and compute the period and preperiod under $f$ 
of the point $x$ in $P_k(S_N)$ such that $x[0,k-1]$ is the chosen
block. Here, by definition   
the period $p$ of $x$ is the eventual period: the smallest 
$j>0$ such that for some $k\geq 0$, $f^k(x)=f^{k+j}(x)$.
So, $p$ is the length of the $f$-cycle 
into which $x$ is mapped by some $f^k$, $k\geq 0$. 
The preperiod of $x$ is the smallest nonnegative $k$ 
such that $f^k(x)$ is $f$-periodic, that is 
$f^k(x)=f^{k+p}(x)$. For a given sample of $m$ 
points from $P_k(S_2)$, the multiplicity $\mu$  of $p$ 
is the number of times the sampled point has the (eventual) period
$p$.  
In the following tables, we restrict to cellular 
automata $f$ on $N=2$ symbols. 
%$p$ denotes the (eventual) $f$-period of a sampled point. 
%$L$  denotes the maximum $f$- period of the  points n
%sampled from  
%In the tables for FProbPeriod, row $k$ corresponds to the 
%points sampled from $P_k(S_2)$. 
% $L$; and $M$ is the maximum preperiod. 
%
%This table  is based at each $k$ on a random sample by 
%FProbPeriod of $m$  blocks of length $k$ on symbols $0,1$. 
%Each block corresponds to a point $x$ in $P_k(S_2)$. 
%It is a table with multiplicities $\mu$ 
%of the (eventual) $B$-periods of the 
%shift-periodic points corresponding to the chosen blocks,   
%where $B$ is the c.a. 
% $B=x_0 + x_1x_2$. It is known that the 
%periodic points of $B$ are dense.

\begin{table}[ht]
\begin{center}
\begin{tabular}{|| r|r|r || r|r|r || r|r|r||}
\hline  
$k$ &$p$ &$\mu$ & $k$ & $p$ &$\mu$ & $k$ & $p$ &$\mu$    \\
\hline 
1 &    1        &     10  &          15&    1        &     1 &  27 &   3402     &     4 \\
2 &    1        &     10  &            &    399      &     3 &	   &   12096    &     3 \\
3 &    1        &     10  &            &    450      &     3 &	   &   218835   &     2 \\
4 &    1        &     4   &            &    530      &     1 &	   &   242352   &     1 \\
  &    4        &     6   &            &    1095     &     2 &	28 &   882      &     1 \\
5 &    1        &     4   &          16 &   120      &     7 &	   &   32144    &     2 \\
  &    15       &     6   &             &   2688     &     3 &	   &   57036    &     7 \\
6 &    1        &     10  &          17 &   2533     &     9 &	29 &   98223    &     2 \\
7 &    1        &     4   &             &   3230     &     1 &	   &   193256   &     1 \\
  &    49       &     6   &          18 &   296      &     1 &	   &   340286    &    3 \\
8 &    1        &     2   &             &   324      &     9 &	   &   504252    &    4 \\
  &    4        &     2   &          19 &   1        &     1 &	30 &   17580     &    1 \\
  &    120      &     6   &             &   13471    &     9 &	   &   161721    &    8 \\
9 &    1        &     2   &          20 &   5240     &     3 &	   &   212670    &    1 \\
  &    9        &     1   &             &   20500    &     1 &	31 &   2228621   &    10\\
  &    54       &     7   &             &   21240    &     6 &	32 &   473792    &    10\\
10&    410      &     10  &          21 &   266      &     10&	33 &   74439     &    7 \\
11&    1        &     1   &          22 &   176      &     2 &	   &   313984    &    3 \\
  &    176      &     9   &             &   14344    &     4 &	34 &   2533   &  10 \\
12&    1        &     8   &             &   32428    &     4 &	35 &   1635074   &6    	\\
  &    56       &     2   &          23 &   1        &     2 &	   &   4485250   &1	\\
13&    1        &     1   &             &   622      &     1 &	   &   14840595  &3	\\
  &    143      &     3   &             &   9108     &     7 &	36 &   152       &2	\\
  &    403      &     2   &          24 &   1        &     2 &	   &   324       &2	\\
  &    416      &     4   &             &   12432    &     4 &	   &   22974     &1	\\
14&    49       &     1   &             &   20256    &     4 &	   &  1700772    &1 \\
  &    448      &     7   &          25 &   61830    &     4 &	   &   4191696	 &4	 \\
  &    882      &     2   &             &   104425   &     6 &	37 &	1365226	 &1	\\
  &		&	  &	    26 &   143      &     1  &	   &	7065594	 &5	\\
  &         	&	  &	       &   6994     &     9  &	   &39209196&4	\\
\hline				      
\end{tabular}
\vskip .1in
\caption{FProbPeriod output for 
the c.a.  $B=x_0 + x_1x_2$ on two symbols, with sample size
 $m=10$. 
%Each block corresponds to a point $x$ in $P_k(S_2)$. 
%The multiplicity $\mu$ is the number of these 10 points 
%which under $B$ eventually fall into a $B$-cycle of 
%the given period. 
% It is known that the 
%periodic points of $B$ are dense.
}\label{fppBcycles}
\end{center}
\end{table}

\begin{table}[ht]
\begin{center}
\begin{tabular}{|| r|r|r|r|r ||  r|r|r|r|r ||  r|r|r|r|r || }
\hline  
$k$ &$p_{30}$ &$\mu_{30}$ & $p_{10}$ &$\mu_{10}$ & 
$k$ &$p_{30}$ &$\mu_{30}$ & $p_{10}$ &$\mu_{10}$ & 
$k$ &$p_{30}$ &$\mu_{30}$ & $p_{10}$ &$\mu_{10}$  \\
\hline 
 1  &   1 &     30  &   1   &  10 & 15 &  1    &	4   &       &   & 23 &   1     &   1  &            &    \\  \cline{1-5}
 2  &   1 &     30  &   1   &  10 &    &  180  &	1   &       &&       &   622   &   2   & 622       & 1   \\	 \cline{1-5}
 3  &   1 &     30  &   1   &  10 &    &  399  &	10  &  399  &6   &    &   9108  &   27 &  9108     & 9    \\	 \cline{1-5}\cline{11-15} 
 4  &   1 &     13  &   1   &  5  &    &  450  &	8   &  450  &1       &24 &  1    &   4  &          &    \\	 
    &	4  &    17  &  4    & 5   &    &  530  &	2   &  	    &    &    &   184   &   1   &          &    \\	 \cline{1-5} 
 5  &   1 &     12  &   1   &  4  &    &  1095 &	5   &1095   &3    &    &   2330  &   1  &   2330   &  1  \\	 \cline{6-10} 
    &	15 &    18  &  15   & 6   & 16 &  1    &	2   &       &    &    &   7440  &   6   &   7440   &  1  \\	 \cline{1-5} 
 6  &   1 &     30  &   1   &  10 &    &  4    &	1   &       &   &    &   12432 &   10 &     12432  &  5   \\	 \cline{1-5} 
 7  &   1 &     9   &   1   &  3  &    &  120  &	12  &	120 &3      &    &   20256 &  8 &   20256  &  3   \\	 \cline{11-15}  
    &	49 &    21  &  49   & 7   &    &  2688 &	15  &   2688&7    & 25 &  4325   &   1  &          &    \\	 \cline{1-5}\cline{6-10}   
 8  &   1 &     2   &   1   &  3  & 17 &  1    &	2   &       &&    &   13015 &  1   &           &    \\	 
    &	4  &    6   &  4    & 2   &    &  2533 &	26  &  2533 &9   &    &   61830 &   14 &  61830    & 4   \\	 
    &	120&    22  &  120  & 5   &    &  3230 &	2   &  3230 &1    &    &  73175  &   4  & 73175    & 1    \\	 \cline{1-5}\cline{6-10} 
 9  &   1 &     3   &       &     & 18 &  1    &	1   &  54   &1   &    &   104425&  10  &  104425   & 5    \\	 \cline{11-15}
    &	54 &    27  &  54   &10   &    &   216 &	 2  &       &   & 26 &    6994 &   30 &   6994        &  10 \\  \cline{1-5}\cline{11-15}
 10 &   1 &     2   &       &     &    &  296  &	5   &	    &    & 27 &    3402 &    9 &    3402       &   3 \\  
    &	15 &    3   &   1   &  2  &    &  324  &	22  &324    &9    &    &   12096 &    9 &    12096      &   3 \\  \cline{6-10}	    
    &	410&    25  &  410  & 8   & 19 &  4161 &	1   &  1    &1   &    &   218835&    10&    218835     &   1 \\  \cline{1-5}
 11 &   1 &     9   &             1   &  2  &    &  4180 &	2   &  4180 &1   &    &   242352&    2 &    242352     &   3 \\	 \cline{11-15}   
    &	176&    21  &  176  & 8   &    & 13471 & 27         & 13471 &8   &28  &  448    &   5  &         &    \\	 \cline{1-5}\cline{6-10}  
 12 &   1 &     27  &   1   &  10 & 20 &  15   &	1   &       &   &    &   882   &   5  &    882   &  4  \\	   
    &	56 &    3   &       &	  &    &  5240 &	5   & 5240  &5   &    &   32144 &   6  &  32144  &  1  \\  \cline{1-5}
 13 &   1 &     1   &       &      &  &  20500&	7           & 20500 &1      &    &   57036 & 14 &  57036 &  5   \\	 \cline{11-15}  
    &	10 &    2   &   10  &  1  &    & 21240 & 17         & 21240 &4&29  &  98223  &  4  &         &    \\	 \cline{6-10}   
    &	143&    13  &  143  & 2   & 21 &  1    &	1   &       &&    &   193256&   7  &         &    \\	 
    &	403&    5   &  403  & 1   &    & 266   & 26         &266    &10  &&   340286& 13   &   340286& 3    \\	 
    &	416&    9   &  416  & 6   &                & 11865 & 3      &	 && & 504252&  6   &   504252&  7   \\  	 \cline{1-5}\cline{6-10}\cline{11-15}
 14 &   49&     2   &   &       & 22 &  176  &	6   &  176  &4   & 30 &   1     &    1 &          &    \\	    
    &	161&    4   &  161  & 3   &    & 1067  & 1          &       &   &    &  450    &   5  &  &  \\	   
    &	448&    11  &  448  & 4   &    & 14344 & 9          &	    &    &    &   1800  &    1 &1800 &1 \\	 
    &	490&    1   &  490  & 1   &    & 32428 & 14         &32428  &6&    &   127995&    2 & 	       &     \\	
    &   882   & 12   & 882  & 2    &    &       &            &&&&132720&1&&\\    
    &	&      &       &	  &    &       &           &    &     &      &   161721&   19 &    161721      &  7           \\  
    &		    &   	  &	&    &       &      &     &    &      &  & 212670&    1 & 212670     &   2      \\  
\hline				      
\end{tabular}
\vskip .1in
\caption{
Table for $B=x_0+x_1x_2$ constructed as for Table \ref{fppCcycles}, 
 for sample sizes 10 and 30 for FProbPeriod. 
The longest orbit length $p$ found is the same for both sample
sizes, except for $k=12$ and $k=21$. 
Exact cycle data for the map $B$, through shift period $k=22$, 
is given in  Table \ref{fpBcycles}. 
Note that for $k\leq 22$, 
only at  $k=12$ did the size-30 
probabilistic sample miss the largest $B$ period.
} 
\label{fppBsamplecycles}
\end{center}
\end{table}

%\clearpage

\begin{table}[ht]
\begin{center}
\begin{tabular}{||r| r|r || r|r ||| r|r || r|r||}
\hline 
map &A & B & $C$     & $E$      & G & H & J & K  \\
\hline 
$k$ &$L^{1/k}$ &$L^{1/k}$ &$L^{1/k}$ &$L^{1/k}$&$L^{1/k}$ &$L^{1/k}$&$L^{1/k}$ &$L^{1/k}$ \\ 
\hline                                                       
1  & 1.00  & 1.00 & 1.00   & 1.00  & 1.00   & 1.00    & 1.00       & 1.00  \\    		 
2  & 1.00  & 1.00 & 1.00   & 1.00  & 1.00   & 1.00    & 1.00       & 1.00  \\  
3  & 1.44  & 1.00 & 1.44   & 1.00  & 1.00   & 1.44    & 1.44       & 1.00  \\      
4  & 1.00  & 1.41 & 1.18   & 1.41  & 1.00   & 1.00    & 1.41       & 1.00  \\      
5  & 1.71  & 1.71 & 1.71   & 1.58  & 1.37   & 1.37    & 1.58       & 1.37  \\      
6  & 1.34  & 1.00 & 1.20   & 1.00  & 1.51   & 1.20    & 1.61       & 1.61  \\  
7  & 1.32  & 1.74 & 1.66   & 1.54  & 1.32   & 1.60    & 1.32       & 1.60  \\  
8  & 1.00  & 1.81 & 1.63   & 1.18  & 1.00   & 1.29    & 1.18       & 1.62  \\  
9  & 1.58  & 1.55 & 1.27   & 1.22  & 1.48   & 1.29    & 1.58       & 1.58  \\  
10 & 1.40  & 1.82 & 1.49   & 1.46  & 1.00   & 1.31    & 1.25       & 1.31  \\   
11 & 1.69  & 1.60 & 1.57   & 1.55  & 1.53   & 1.43    & 1.46       & 1.57  \\  
12 & 1.23  & 1.39 & 1.16   & 1.00  & 1.30   & 1.34    & 1.30       & 1.36  \\  
13 & 1.67  & 1.59 & 1.63   & 1.49  & 1.53   & 1.48    & 1.53       & 1.66  \\  
14 & 1.20  & 1.62 & 1.38   & 1.35  & 1.37   & 1.38    & 1.39       & 1.26  \\  
15 & 1.19  & 1.59 & 1.50   & 1.39  & 1.40   & 1.42    & 1.44       & 1.48  \\  
16 & 1.00  & 1.63 & 1.50   & 1.47  & 1.40   & 1.13    & 1.54       & 1.63  \\  
17 & 1.38  & 1.60 & 1.68   & 1.55  & 1.45   & 1.52    & 1.52       & 1.55  \\  
18 & 1.30  & 1.37 & 1.58   & 1.41  & 1.33   & 1.40    & 1.53       & 1.59  \\  
19 & 1.62  & 1.64 & 1.60   & 1.45  & 1.50   & 1.51    & 1.47       & 1.49  \\  
20 & 1.22  & 1.64 & 1.54   & 1.50  & 1.32   & 1.18    & 1.12       & 1.56  \\  
21 & 1.21  & 1.30 & 1.48   & 1.36  & 1.15   & 1.29    & 1.48       & 1.52  \\  
22 & 1.34  & 1.60 & 1.51   & 1.46  & 1.23   & 1.49    & 1.40       & 1.35  \\  
23 & 1.39  & 1.48 & 1.57   & 1.53  & 1.44   & 1.55    & 1.53       & 1.58  \\  
24 & 1.14  & 1.51 & 1.54   & 1.45  & 1.39   & 1.14    & 1.46       & 1.46  \\  
25 & 1.50  & 1.58 & 1.48   & 1.50  & 1.23   & 1.53    & 1.44       & 1.56  \\  
26 & 1.32  & 1.40 & 1.56   & 1.46  & 1.33   & 1.44    & 1.43       & 1.29  \\  
27 & 1.42  & 1.58 & 1.54   & 1.50  & 1.43   & 1.46    & 1.28       & 1.58  \\  
28 & 1.12  & 1.47 & 1.43   & 1.47  & 1.30   & 1.31    & 1.48       & 1.44  \\  
29 & 1.56  & 1.57 & 1.54   & 1.49  & 1.31   & 1.35    & 1.51       & 1.63  \\  
30 & 1.12  & 1.50 & 1.59   & 1.44  & 1.33   & 1.37    & 1.45       & 1.64  \\  
31 & 1.11  & 1.60 & 1.59   & 1.41  & 1.47   & 1.32    & 1.49       & 1.63  \\  
32 & 1.00  & 1.50 & 1.60   & 1.26  & 1.31   & 1.46    & 1.48       & 1.54  \\  
33 & 1.23  & 1.46 & 1.63   & 1.44  & 1.37   & 1.43    & 1.49       & 1.54  \\  
34 & 1.20  & 1.25 & 1.45   & 1.48  & 1.35   & 1.48    & 1.48       & 1.56  \\  
35 & 1.26  & 1.60 & 1.44   & 1.48  & 1.42   & 1.43    & 1.48       & 1.59  \\  
36 & 1.16  & 1.52 & 1.50   & 1.41  & 1.27   & 1.29    & 1.38       & 1.46  \\  
37 & 1.49  & 1.60 & 1.55   & 1.46  & 1.30   & 1.44    & 1.45       & 1.57  \\  
\hline
\end{tabular}
\vskip .1in
\caption{This table  is based at each $k$ on a random sample 
of 10 blocks of length $k$ on two symbols. $L$ is the 
maximum period  from the sample.  \newline 
 $A=x_0+x_1$, linear, and    
$B=x_0 +x_1x_2$, permutative.   \newline 
$C= B\circ B_{rev}$, nonclosing, where 
$B_{rev}=x_0x_1+x_2$.  \newline 
$G= x_{-1}+x_0x_1+x_2$,   bipermutative, nonlinear. \newline 
$E=B\circ A$,  degree 2, left closing, not 
right closing.  \newline 
$J=A\circ U$, where 
 $U= x_0 + x_{-2}(1+x_{-1})x_1x_2$ is invertible.   \newline 
%where  $U$ is invertible. which applies  the flip to the symbol 
%in the $*$ space of the  frame $1 0 * 1 1$.    which equals 
%$ x_0 + x_{-2}x_1x_2 +  x_{-2}x_{-1}x_1x_2$. 
$H=A\circ A \circ U$, and   
$K=B\circ U$.  
}\label{fpptable}
\end{center}
\end{table}

\begin{table}[ht]
\begin{center}
\begin{tabular}{||r| r|r || r|r ||| r|r || r|r||}
\hline 
map &A & B & C     &  E     & G & H & J & K \\
\hline 
$k$ &$M^{1/k}$ &$M^{1/k}$ &$M^{1/k}$ &$M^{1/k}$&$M^{1/k}$ &$M^{1/k}$&$M^{1/k}$ &$M^{1/k}$ \\ 
\hline 
1  & 1.00   & 1.00   & 1.00   & 1.00   & 0.00  & 1.00  & 1.00 & 1.00     \\ 
2  & 1.41   & 1.00   & 1.00   & 1.00   & 1.00  & 1.00  & 1.41 & 1.00     \\  
3  & 1.00   & 1.25   & 1.00   & 1.25   & 1.25  & 1.00  & 1.00 & 1.25     \\      
4  & 1.41   & 1.00   & 1.00   & 1.18   & 1.31  & 1.18  & 1.31 & 1.00     \\      
5  & 1.00   & 1.00   & 1.00   & 1.00   & 1.31  & 1.24  & 1.14 & 1.00     \\      
6  & 1.12   & 1.41   & 1.20   & 1.25   & 1.20  & 1.12  & 1.12 & 1.25     \\  
7  & 1.00   & 1.21   & 1.16   & 1.10   & 1.34  & 1.16  & 1.32 & 1.25     \\  
8  & 1.29   & 1.09   & 1.00   & 1.33   & 1.36  & 1.18  & 1.36 & 1.31     \\  
9  & 1.00   & 1.25   & 1.25   & 1.27   & 1.19  & 1.24  & 1.24 & 1.24     \\  
10 & 1.07   & 1.14   & 1.32   & 1.11   & 1.30  & 1.21  & 1.30 & 1.30     \\   
11 & 1.00   & 1.44   & 1.23   & 1.26   & 1.28  & 1.25  & 1.25 & 1.42     \\  
12 & 1.12   & 1.44   & 1.37   & 1.34   & 1.26  & 1.20  & 1.32 & 1.38     \\  
13 & 1.00   & 1.34   & 1.27   & 1.32   & 1.21  & 1.25  & 1.23 & 1.41     \\  
14 & 1.05   & 1.41   & 1.37   & 1.32   & 1.30  & 1.25  & 1.35 & 1.42     \\  
15 & 1.00   & 1.39   & 1.37   & 1.30   & 1.26  & 1.34  & 1.32 & 1.41     \\  
16 & 1.18   & 1.46   & 1.39   & 1.34   & 1.28  & 1.25  & 1.29 & 1.32     \\  
17 & 1.00   & 1.47   & 1.34   & 1.33   & 1.28  & 1.26  & 1.26 & 1.41     \\  
18 & 1.03   & 1.46   & 1.34   & 1.31   & 1.23  & 1.22  & 1.30 & 1.40     \\  
19 & 1.00   & 1.44   & 1.34   & 1.32   & 1.25  & 1.34  & 1.32 & 1.42     \\  
20 & 1.07   & 1.39   & 1.37   & 1.28   & 1.28  & 1.27  & 1.36 & 1.37     \\  
21 & 1.00   & 1.46   & 1.41   & 1.36   & 1.33  & 1.35  & 1.33 & 1.48     \\  
22 & 1.03   & 1.44   & 1.39   & 1.34   & 1.26  & 1.23  & 1.31 & 1.48     \\  
23 & 1.00   & 1.47   & 1.40   & 1.33   & 1.29  & 1.34  & 1.31 & 1.44     \\  
24 & 1.09   & 1.43   & 1.43   & 1.35   & 1.32  & 1.28  & 1.29 & 1.44     \\  
25 & 1.00   & 1.44   & 1.41   & 1.33   & 1.30  & 1.27  & 1.30 & 1.47     \\  
26 & 1.02   & 1.45   & 1.43   & 1.31   & 1.28  & 1.26  & 1.31 & 1.45     \\  
27 & 1.00   & 1.45   & 1.43   & 1.33   & 1.31  & 1.37  & 1.37 & 1.45     \\  
28 & 1.05   & 1.46   & 1.42   & 1.32   & 1.30  & 1.29  & 1.34 & 1.46     \\  
29 & 1.00   & 1.42   & 1.41   & 1.31   & 1.31  & 1.33  & 1.32 & 1.44     \\  
30 & 1.02   & 1.41   & 1.44   & 1.34   & 1.31  & 1.28  & 1.35 & 1.47     \\  
31 & 1.00   & 1.44   & 1.43   & 1.38   & 1.30  & 1.36  & 1.34 & 1.45     \\  
32 & 1.11   & 1.46   & 1.43   & 1.37   & 1.30  & 1.31  & 1.34 & 1.45     \\  
33 & 1.00   & 1.46   & 1.42   & 1.35   & 1.27  & 1.34  & 1.33 & 1.46     \\  
34 & 1.03   & 1.47   & 1.24   & 1.34   & 1.33  & 1.32  & 1.34 & 1.47     \\  
35 & 1.00   & 1.45   & 1.56   & 1.37   & 1.31  & 1.34  & 1.37 & 1.46      \\  
36 & 1.04   & 1.44   & 1.46   & 1.36   & 1.29  & 1.32  & 1.36 & 1.47      \\  
37 & 1.01   & 1.47   & 1.46   & 1.35   & 1.32  & 1.34  & 1.35 & 1.46      \\  
\hline
\end{tabular}
\vskip .1in
\caption{This table  is based at each $k$ on a random sample 
of 10 blocks of length $k$ on two symbols. $M$ 
is the maximum preperiod from the sample.  The 
maps are the same as used in Table \ref{fpptable}
\newline 
(Remark: For $x\in P_k(S_2)$
with $k=q2^j$ with $q$ odd, the point $A^{j+1}(x) $
must be periodic under $f$, so $M\leq j+1$.) 
}\label{fppptable}
\end{center}
\end{table}

\begin{table}[ht]
\begin{center}
\begin{tabular}{|| r|r|r || r|r|r || r|r|r  || r|r|r || }
\hline  
$k$ &$p$ &$\mu$ & $k$ & $p$ &$\mu$ & $k$ & $p$
&$\mu$ & $k$ & $p$ &$\mu$       \\
\hline 
1   & 1      & 10  & 12     &1     & 10 &20    &3      & 1   &30  & 31     &1   \\
2   & 1      & 10  & 13     &47    & 6  &      &2790   & 9   &    & 82531  &9    \\
3   & 1      & 10  &        &52    & 4  &21    &573    & 10  &31  & 57747  &10  \\
4   & 1      & 10  & 14     &1     & 1  &22    &11     & 3   &32  & 85     &2   \\
5   & 1      & 1   &        &5     & 3  &      &519    & 3   &    & 91     &2   \\
    & 3      & 9   &        &13    & 3  &      &9658   & 4   &    & 234649 &6   \\
6   & 1      & 5   &        &47    & 2  &23    &1499   & 8   &33  & 3452570&6   \\
    & 3      & 5   &        &49    & 1  &      &9384   & 2   &    &10357710&4   \\
7   & 1      & 3   & 15     &31    & 7  &24    &1      & 1   &34  &  1717  &7   \\
    & 5      & 7   &        &145   & 3  &      &35160  & 9   &    &10574   &1   \\
8   & 1      & 5   & 16     &29    & 1  &25    &20475  & 10  &    &999056  &2  \\
    & 13     & 5   &        &85    & 9  &26    &441    & 1   &35  &572068  &1  \\
9   & 1      & 1   & 17     &101   & 3  &      &9401   & 9   &    &2860340 & 1  \\
    & 9      & 9   &        &399   & 7  &27    &4543   & 5   &    &3262280 & 8  \\
10  & 1      & 1   & 18     &1     & 1  &      &19710  & 1   &36  &56      & 2  \\
    & 3      & 7   &        &455   & 9  &      &113643 & 4   &    &4095    & 1  \\
    & 5      & 1   & 19     &1     & 1  &28    &5      & 1   &    &729537  & 3   \\
    & 11     & 1   &        &401   & 4  &      &1260   & 9   &    &908910  & 1   \\
11  & 11     & 6   &        &2755  & 3  &29    &277298 & 10  &    &2188611 & 3   \\
    & 143    & 4   &        &7125  & 2  &      &       &     &37  &5881335 & 2  \\
    &        &     &        &      &    &      &       &     &    &12081277& 8  \\
\hline				      
\end{tabular}
\vskip .1in
\caption{
This table  is based at each $k$ on a random sample by 
FProbPeriod 
of 10 blocks of length $k$ on symbols $0,1$. It is a table of 
the resulting periods $p$  with multiplicities $\mu$
for  the nonclosing  map $C=B\circ B_{rev}$ where $B=x_0 + x_1x_2$
and $B_{rev}=x_0x_1+x_2$. Here, a sampled block $x[0,k-1]$ determines a 
point $x$ of period $k$,  and the period $p$ is by 
definition the eventual period of $x$ under iteration by $C$.
The multiplicity $\mu$ for a given period is the number of samples 
for which $x$ has that period.}\label{fppCcycles}
\end{center}
\end{table}

\begin{table}[ht]
\begin{center}
\begin{tabular}{|| r|r|r || r|r|r || r|r|r  || r|r|r || r|r|r ||}
\hline  
$k$ &$p$ &$\mu$ & $k$ & $p$ &$\mu$ & $k$ & $p$
&$\mu$ & $k$ & $p$ &$\mu$   & $k$ & $p$ &$\mu$    \\
\hline 
1   &  1    &  1      &  11   & 11   & 7   & 19  &2755   &3   & 27 & 113643 &6   &34&1717 & 10  \\
2   &  1    &  1      &       & 143  & 3   &     &7125 	 &5   &    & 122661 &4   &35  &3262280 &8  \\
3   &  1    &  7      &  12   & 1    & 9   &     &7619 	 &2   & 28 & 35     &2   &  &6886355 & 2 \\
    &  3    &  3      &       & 6    & 1   & 20  &1395 	 &9   &    & 180    &5   &36&56 & 2 \\
4   &  1    &  6      &  13   & 13   & 2   &     &5780 	 &1   &    & 1260   &2   &  &504 & 2 \\
    &  2    &  4      &       & 52   & 5   & 21  &4011 	 &10  &    & 26124  &1   &  &729537 & 3 \\
5   &  1    &  2      &       & 611  & 3   & 22  &11   	 &1   & 29 & 277298 &10  &  &2188611 & 3 \\
    &  5    &  3      &  14   & 7    & 2   &     &878  	 &5   & 30 & 5205   &1   &37&3768043  & 1 \\
    &  15   &  5      &       & 35   & 2   &     &5709 	 &2   &    & 137190 &1   &  &5881335 & 2 \\
6   &  1    &  7      &       & 91   & 6   &     &9658 	 &2   &    & 1237965&8   &  &12081277 & 7 \\
    &  3    &  3      &  15   & 87   & 2   & 23  &9384 	 &2   & 31 & 457777 &1   &  & &  \\
7   &  7    &  2      &       & 465  & 8   &     &34477	 &8   &    & 1790157&9   &  & &  \\
    &  35   &  8      &  16   & 680  & 7   & 24  &1      &2   & 32 & 680    &1   &  & &   \\
8   &  1    &  1      &       & 728  & 3   &     &11720  &5   &    & 728    &3   &  & &   \\
    &  4    &  2      &  17   & 1717 & 5   &     &35160  &3   &    & 267824 &2   &  & &   \\
    &  52   &  7      &       & 6783 & 5   & 25  &4095   &8   &    & 3754384&4   &  & &   \\
9   &  1    &  1      &  18   & 9    & 1   &     &20475  &2   & 33 & 3452570&5    &  & &    \\
    &  9    &  9      &       & 56   & 1   & 26  &52     &2   &    &10357710&5    &  & &    \\
10  &  15   &  7      &       & 4095 & 8   &     &122213 &8   &    &        &    &  & &    \\
    &  55   &  3      &       & &    &     &    	 &    &    &        &    &  & &    \\
\hline				      
\end{tabular}
\vskip .1in
\caption{
This table  is constructed just as Table \ref{fppCcycles} 
was, except for the following: the data is for the map $D$ 
which is the
map $C$ composed with $(S_2)^{-2}$, i.e., 
$D$ is the composition of $x_0 +x_1x_2$ with   $x_{-2}x{-1}+x_0$.}
%The cycle lengths look little more ``random'' for $D$.Maybe??} 
\label{fppDcycles}
\end{center}
\end{table}

\begin{table}[ht]
\begin{center}
\begin{tabular}{|| r|r|r|r || r|r|r|r || r|r|r|r ||}
\hline  
$k$ &$p$ &$\mu_{30}$ & $\mu_{10}$ & 
$k$ &$p$ &$\mu_{30}$ & $\mu_{10}$ & 
$k$ &$p$ &$\mu_{30}$ & $\mu_{10}$   \\
\hline 
1   & 1      & 30       & 10    &17  & 255      &30   &10 & 36 & 252     &30&10\\  
2   & 1      & 30       & 10    &18  & 126      &30   &10 & 37 & 3233097 &30&10\\
3   & 1      & 7        & 3     &19  & 9709     &30   &10 & 38 & 19418   &30&10\\
    & 3      & 23       & 7     &20  & 30       &1    &1  & 39 & 4095    &30&10\\
4   & 1      & 30       & 10    &    & 60       &29   &9  & 40 & 120     &30&10\\
5   & 15     & 30       & 10    &21  & 63       &30   &10 & 41 & 41943   &30&10\\
6   & 1      & 3        &  0    &22  & 682      &30   &10 & 42 & 126     &30&10\\
    & 3      & 9        & 3     &23  & 2047     &30   &10 & 43 & 5461    &30&10\\
    & 6      & 18       & 7     &24  & 24       &30   &10 & 44 & 1364    &30&10\\
7   & 7      & 30       & 1     &25  & 25575    &30   &10 & 45 & 4095    &30&10\\
8   & 1      & 30       & 9     &26  & 1638     &30   &10 & 46 & 4094    &30&10\\
9   & 63     & 30       & 10    &27  & 13797    &30   &10 & 47 & 8388607 &30&10\\
10  & 15     & 1        & 10    &28  & 28       &30   &10 & 48 & 48      &30&10\\
    & 30     & 29       & 10    &29  & 475107   &30   &10 & 49 & 2097151 &30&10\\
11  & 341    & 30       & 10    &30  & 30       &30   &10 &  &  &&  \\
12  & 12     & 30       & 10    &31  & 31       &30   &10 &  &  &&  \\
13  & 819    & 30       & 10    &32  & 1        &30   &10 &  &  &&  \\
14  & 14     & 30       & 10    &33  & 1023     &30   &10 &  &  &&  \\
15  & 15     & 30       & 10    &34  & 510      &30   &10 &  &  &&  \\
16  & 1      & 30       & 10    &35  & 4095     &30   &10 &  &  &&  \\
\hline				      
\end{tabular}
\vskip .1in
\caption{
Table for $A=x_0+x_1$ constructed as for Table \ref{fppCcycles}, 
 for sample sizes 10 and 30 for FProbPeriod. Both sample sizes work 
well  for this linear map. 
} 
\label{fppAsamplecycles}
\end{center}
\end{table}

%for the nonclosing map $C= B_{rev}\circ B$ 
%for $N=2$ symbols for sample size $m=10$, 
%where $B=x_0 + x_1x_2$ and $B_rev = x_0x_1+x_2$. }\label{fppcycles}
%
%
%
%
%34 &1717      &  10
%35 &3262280  & 8   
%   &    6886355  &2
%36 &56           &2
%   &    504      &2
%   &    729537   &3
%   &    2188611  &3

\begin{table}[ht]
\begin{center}
\begin{tabular}{|| r|r || r|r || r|r || r|r || r|r || r|r || r|r ||}
\hline  
k  &   Pre.        &k    &   Pre.  &k    & Pre.        &k  &Pre. &k&Pre.&k&Pre.  \\ 
\hline 
18 &   0    	   &	21 &   216 & 	24 &     394 & 27  &   2669&30   & 7988    &33  &   	100664 \\   
   &   45   	   &	   &   232 & 	   &   1032  &     &   3329  & 	   &   8537    &  & 	108279 \\	     
   &   170  	   &	   &   254 & 	   &   1627  &     &   5641  & 	   &   10545   &  & 	146726 \\	     
   &   226  	   &	   &   578 & 	   &   1834  &     &   9244  & 	   &   12479   &  & 	149017 \\	     
   &   293  	   &	   &   595 & 	   &   2880  &     &   9697  & 	   &   13204   &  & 	157529 \\	     
   &   362  	   &	   &   756 & 	   &   3145  &     &   11036 & 	   &   13676   &  & 	161009 \\	     
   &   506  	   &	   &   1174& 	   &   3396  &     &   11583 & 	   &   16210   &  & 	188071  \\	     
   &   556  	   &	   &   1313& 	   &   3905  &     &   13921 & 	   &   16315   &  & 	196758 \\	     
   &   751  	   &	   &   3058& 	   &   5215  &     &   14745 & 	   &   24373   &  & 	240207\\
   &   102  	   &	   &   3186& 	   &   5517  &     &   23357 & 	   &   31106   &  & 	270439 \\	     
\hline		    								           
19 &   0    	   &	22 &   73  & 	25 &     57  &  28 &   6000  & 	31 &   4184    &34&74577 	\\	     
   &   9    	   &	   &   350 & 	   &   404   &     &   7790  & 	   &   11132   &  &80674 	\\	     
   &   145  	   &	   &   381 & 	   &   463   &     &   14999 & 	   &   16211   &  &161429	\\	     
   &   201  	   &	   &   522 & 	   &   750   &     &   18569 & 	   &   30108   &  &193935	\\	     
   &   291  	   &	   &   587 & 	   &   1361  &     &   23067 & 	   &   30661   &  &209842	\\	     
   &   301  	   &	   &   1163& 	   &   1547  &     &   25108 & 	   &   32789   &  &283852	\\	     
   &   490  	   &	   &   1480& 	   &   1671  &     &   28943 & 	   &   39001   &  &360139	\\	     
   &   547  	   &	   &   1788& 	   &   2491  &     &   31291 & 	   &   57399   &  &400913	\\	     
   &   658  	   &	   &   3131& 	   &   9415  &     &   37745 & 	   &   70339   &  &452695	\\
   &   105  	   &	   &   3247& 	   &   9531  &     &   47285 &    &   99228    &  &521740	\\	     
\hline		    								           
20 &   0    	   &     23&      0&     26&     428 &   29&    2702 &    32&	8955   &35&38461 	\\	     
   &   44   	   &	   &      0&       &    7866 &     &    9704 &      &  14007   &  &51383 	\\		     
   &   181  	   &	   &   671 &       &    9538 &     &    10239&      &	19792  &  &71921 	\\		     
   &   230         &	   &   2543&       &   11744&      &   10783 &      &	48308  &  &82039 	\\		     
   &   257  	   &	   &   2791&       &    12056&     &    11194&      &	49195  &  &122647	\\		     
   &   290  	   &	   &   2835&       &    13070&     &    13472&      &	65168  &  &261207	\\		     
   &   445  	   &	   &   2900&       &    15762&     &    19429&      &	65824  &  &416233	\\		     
   &   556  	   &	   &   3123&       &    16164&     &    22208&      &	72902  &  &492662	\\		     
   &   622  	   &	   &   6492&       &    16999&     &    31392&      &	126309 &  &529999	\\		     
   &   755  	   &	   &   7854&       &    17346&     &    31856&      &   184781 &  &555899	\\ 
\hline
\end{tabular}
\vskip .1in
\caption{This table  is based at each $k$ on a random sample 
of 10 blocks of length $k$ on symbols $0,1$. It is a table of 
the preperiods (Pre.) seen by FProbPeriod in a sample for the 
 map $B=x_0 +x_1x_2$. 
}\label{preperB}4
\end{center}
\end{table}

\begin{table}[ht]
\begin{center}
\begin{tabular}{|| r|r || r|r || r|r || r|r || r|r || r|r || r|r ||}
\hline  
k  &   Pre.        &k    &   Pre.  &k    & Pre.        &k  &Pre. &k&Pre.&k&Pre.  \\ 
\hline 
18 &   0      	   &	21 & 134    & 	24 & 316    & 27  & 2016    &30    &1862      &33&1790     	 \\   
   &   15     	   &	   & 159    & 	   & 969    &     & 2046    & 	   &2020      &  &6742   	 \\	     
   &   17     	   &	   & 286    & 	   & 1012   &     & 2543    & 	   &2159      &  &8823   	 \\	     
   &   26     	   &	   & 305    & 	   & 1027   &     & 3793    & 	   &6736      &  &19037  	 \\	     
   &   53     	   &	   & 551    & 	   & 1541   &     & 5981    & 	   &15523     &  &19158  	 \\	     
   &   98     	   &	   & 600    & 	   & 1954   &     & 10577   & 	   &16289     &  &19234  	 \\	     
   &   113    	   &	   & 696    & 	   & 2800   &     & 13024   & 	   &29590     &  &33862  	  \\	     
   &   177    	   &	   & 748    & 	   & 3054   &     & 14611   & 	   &32253     &  &56053  	 \\	     
   &   184    	   &	   & 1379   & 	   & 3371   &     & 15185   & 	   &56898     &  &105560 	\\
   &   198    	   &	   & 1550   & 	   & 5430   &     & 17428   & 	   &59246     &  &119315 	 \\	     
\hline 	    		       		       		     	       
19 &   15     	   &	22 & 7      & 	25 & 1215   &  28 & 982     & 	31 &2176      &34&28492 	\\	     
   &   35     	   &	   & 14     & 	   & 1379   &     & 5042    & 	   &3552      &  &176703	\\	     
   &   48     	   &	   & 202    & 	   & 3305   &     & 5438    & 	   &3716      &  &184042	\\	     
   &   49     	   &	   & 241    & 	   & 3429   &     & 5783    & 	   &16352     &  &203938	\\	     
   &   108    	   &	   & 389    & 	   & 3486   &     & 5890    & 	   &30302     &  &204854	\\	     
   &   119    	   &	   & 436    & 	   & 3673   &     & 8235    & 	   &53475     &  &261574	\\	     
   &   211    	   &	   & 481    & 	   & 3865   &     & 9723    & 	   &55921     &  &282167	\\	     
   &   264    	   &	   & 1003   & 	   & 4716   &     & 10118   & 	   &58664     &  &295105	\\	     
   &   271    	   &	   & 1119   & 	   & 4719   &     & 15474   & 	   &60608     &  &374994	\\
   &   275    	   &	   & 1477   & 	   & 6056   &     & 19004   &      &65768     &  &384473	\\	     
\hline 	    		       				 	       
20 &  7      	   &     23& 177    &    26& 289    &   29& 7160    &    32&2686      &35&15491   	  \\	     
   &  20     	   &	   & 198    &      & 327    &     & 8519    &      &4683      &  &68409   	  \\		     
   &  136    	   &	   & 341    &      & 1127   &     & 9150    &      &9120      &  &76254   	  \\		     
   &  138          &	   & 366    &      & 1301   &     & 16361   &      &9808      &  &134451  	  \\		     
   &  273    	   &	   & 793    &      & 4720   &     & 16679   &      &18985     &  &168960  	  \\		     
   &  280    	   &	   & 1036   &      & 5146   &     & 21158   &      &19562     &  &205781  	  \\		     
   &  422    	   &	   & 1043   &      & 5817   &     & 22279   &      &32081     &  &223889  	  \\		     
   &  474    	   &	   & 1402   &      & 8610   &     & 24581   &      &76390     &  &247831  	  \\		     
   &  478    	   &	   & 2389   &      & 11410  &     & 25450   &      &89515     &  &351861  	  \\		     
   &  614    	   &	   & 2613   &      & 12155  &     & 25558   &      &98677     &  &385646  	  \\ 
\hline
\end{tabular}
\vskip .1in
\caption{This table  is based at each $k$ on a random sample 
of 10 blocks of length $k$ on symbols $0,1$. It is a table of 
the preperiods (Pre.) seen by FProbPeriod in a sample for the 
nonclosing  map $C=B\circ B_{rev}$, where $B=x_0 +x_1x_2$ and 
$B_{rev}=x_0x_1+x_2$.  
}\label{preperC}
\end{center}
\end{table}

\end{document}